\documentclass[preprint,1p,number,sort&compress]{elsarticle}
\usepackage{microtype}
\usepackage{amsmath,amsfonts}
\usepackage{graphicx}
\usepackage{pgfplots}
\usepackage{upquote}
\newcommand{\rmi}{\mathrm{i}}
\newcommand{\rme}{\mathrm{e}}

\newcommand{\CC}{\mathbb{C}}
\newcommand{\bu}{\boldsymbol u}

\newcommand{\bU}{\boldsymbol U}
\newcommand{\bv}{\boldsymbol v}
\newcommand{\bw}{\boldsymbol w}
\newcommand{\bx}{\boldsymbol x}
\newcommand{\bg}{\boldsymbol g}
\newcommand{\bof}{\boldsymbol f}
\DeclareMathOperator{\at2}{atan2}
\DeclareMathOperator{\sech}{sech}
\pgfplotsset{compat=newest}
\usetikzlibrary{plotmarks}
\usetikzlibrary{arrows.meta}
\usepgfplotslibrary{patchplots}
\usepackage{grffile}
\title{Efficient simulation of complex Ginzburg--Landau equations
using high-order exponential-type methods}
\author[1]{Marco Caliari\corref{cor1}}
\ead{marco.caliari@univr.it}

\author[1]{Fabio Cassini}
\ead{fabio.cassini@univr.it}

\cortext[cor1]{Corresponding author}

\affiliation[1]{organization={Department of Computer Science,
    University of Verona},addressline={Strada Le Grazie, 15},
  postcode={37134},
  city={Verona},
country={Italy}}

\usepackage{hyperref}
\begin{document}
\myfooter[L]{}

\begin{abstract}
  In this paper, we consider the task of efficiently computing the numerical
  solution of evolutionary complex Ginzburg--Landau equations
  {\color{black}on Cartesian product domains with homogeneous Dirichlet/Neumann
  or periodic boundary conditions}. To this aim,
  we employ for the time integration
  high-order exponential methods of splitting and
  Lawson type
  {\color{black}with constant time step size}.
  These schemes enjoy favorable stability
  properties and, in particular, do not show restrictions on the time step size
  due to the underlying stiffness of the models. The needed actions of matrix
  exponentials are efficiently realized
by using a tensor-oriented approach that suitably employs the so-called 
  $\mu$-mode product (when the semidiscretization in space is performed 
  with finite differences) or
  with pointwise operations in
  Fourier space (when the model is considered with periodic boundary
  conditions). The overall effectiveness of the approach
  is demonstrated by running simulations on a variety of two- and
  three-dimensional
  (systems of) complex Ginzburg--Landau equations with cubic or
  cubic-quintic nonlinearities,
  which are widely considered in literature to model relevant physical
  phenomena.
{\color{black}
  In fact, we show that high-order exponential-type schemes may
  outperform standard techniques to integrate in time the models under 
  consideration, i.e., the well-known second-order
  split-step method and the explicit
  fourth-order 
  Runge--Kutta integrator, for stringent accuracies.}
\end{abstract}


\maketitle

\section{Introduction}
The evolutionary Complex Ginzburg--Landau (CGL) equation is an ubiquitous model
for many physical scenarios. For instance, such a Partial Differential Equation (PDE)
has applications in nonlinear fiber optics, fluids, and Bose--Einstein 
condensation, especially for the study of the dynamics of the so-called 
dissipative solitons~\cite{MMLLM08} (see also the comprehensive
reviews~\cite{AK02,M22} and references therein).
In this work, we consider the CGL equation in dimensionless units written in
the form
\begin{equation}\label{eq:cqCGL}
  \partial_tu=(\alpha_1+\rmi\beta_1)\Delta u+
  \alpha_2 u+
  (\alpha_3+\rmi\beta_3)\lvert u\rvert^2u+
    (\alpha_4+\rmi\beta_4)\lvert u\rvert^4u,
\end{equation}
where $u\colon[0,T]\times\Omega\to\CC$ represent{\color{black}s}
the complex-valued unknown
function. The physical
meaning of $u$ depends on the specific example under consideration.
The model is completed with an appropriate initial
{\color{black}solution and boundary conditions of homogeneous Dirichlet/Neumann
  or periodic type. Results of existence and uniqueness of solution for
  equation~\eqref{eq:cqCGL} can be found, for instance,
  in References~\cite{D94,DM00,SYY16}.}
The spatial domain $\Omega$, possibly after a truncation,
is assumed to be the Cartesian product of one-dimensional finite intervals,
that is
{\color{black}
$\Omega=(a_1,b_1)\times \cdots \times (a_d,b_d)$}. The parameters
$\alpha_j,\beta_j$, with $j=1,\ldots,4$, are real numbers.
In the following, we
assume that $\alpha_1$ is positive and that $\beta_1$, $\alpha_2$, $\alpha_3$,
and $\beta_3$ are different from zero. If $\alpha_4 = \beta_4 = 0$, the
equation is usually referred to as \emph{cubic complex Ginzburg--Landau} equation,
and it was first introduced in Reference~\cite{SS71} in the context of linearized
perturbations. Whenever $\alpha_4$ and/or $\beta_4$ are different
from zero, the equation is known as
\emph{cubic-quintic complex Ginzburg--Landau}. This type of nonlinearity
was first considered in Reference~\cite{PS84} for the study of
certain soliton dynamics in chemical reactions.
Also, notice that model~\eqref{eq:cqCGL} is in fact a 
generalization of the well-known nonlinear Schr\"odinger equation
that we do
not consider in this manuscript.

A suitable semidiscretization in space of equation~\eqref{eq:cqCGL} by the
method of lines 
on a spatial grid of size $n_1\times\cdots \times n_d$ gives
a system of Ordinary Differential Equations (ODEs) in the form
\begin{subequations}\label{eq:ODEs}
\begin{equation}
\bu'(t) = K\bu(t)+\bg(\bu(t))=\bof(\bu(t)).
\end{equation}
Here $K$ is a complex matrix of size $N\times N$, with
$N=n_1\cdots n_d$, that discretizes the linear operator
$((\alpha_1+\rmi\beta_1)\Delta + \alpha_2)$, while $\bg$ represents the discretization
of the nonlinear part of the PDE. We also assume that the matrix $K$ can be written
as a Kronecker sum, that is
\begin{equation}
K=A_d\oplus A_{d-1} \oplus \cdots \oplus A_1=\sum_{\mu=1}^d A_{\otimes \mu},
\end{equation}
where
\begin{equation}
  A_{\otimes \mu}=I_d\otimes \cdots \otimes I_{\mu+1}\otimes A_\mu
  \otimes I_{\mu-1} \otimes \cdots \otimes I_1.
\end{equation}
\end{subequations}
The matrices $I_\mu$ and $A_\mu$, for $\mu=1,\ldots,d$, are of
size $n_\mu\times n_\mu$, and represent the identity and the discretization
of the one-dimensional linear differential operator along
the $\mu$th direction, that is
\begin{equation*}
A_\mu\approx (\alpha_1+\rmi\beta_1)\partial_{x_\mu x_\mu}+\frac{\alpha_2}{d},
\end{equation*}
respectively. The matrices $A_\mu$ usually incorporate boundary conditions
and they can be sparse or dense,
depending on the specific discretization technique.
Typical examples of discretizations that lead to such a Kronecker structure
are finite differences and spectral methods. 
The former, in particular, is widely employed in the literature when dealing
with non-periodic boundary conditions (see, for instance, 
References~\cite{EHGP01,MMLLM08,WZ13}). On the other hand, if 
equation~\eqref{eq:cqCGL} is coupled with periodic boundary condition, 
the most used space discretization is the Fourier pseudospectral method, which
gives rise to a system of ODEs in the form~\eqref{eq:ODEs} for the Fourier coefficients of the expansion
(we refer to References~\cite{HFDW06,ZY15b,DK18}, among the others, for more
details).
In any case, remark that the resulting system of ODEs is \emph{stiff}, due to the presence
of the second-order spatial differential operator in equation~\eqref{eq:cqCGL}.

Concerning the time discretization, fully explicit methods (usually of 
Runge--Kutta type) are widely employed in the literature,
even if they must take into
account a time step size restriction to guarantee stability
(see, for instance, Reference~\cite{EHGP01}).
However, probably the most popular techniques for the time marching are
based on splitting
schemes, giving rise to Time Splitting Fourier Pseudospectral (TSFP)
methods  or Time Splitting Finite Difference (TSFD) schemes 
(see Reference~\cite{WZ13}, and Reference~\cite{BC13} for an
exhaustive treatment in the field of Bose--Einstein condensates).
In this context, the flow of the original equation is split
into a linear subflow and a nonlinear one.
The linear flow can be approximated
by suitable implicit methods 
when dealing with finite difference discretizations
(see again References \cite{BC13,WZ13}) or solved exactly in Fourier space
(by computing a diagonal matrix exponential).
The nonlinear flow can typically be solved analytically, otherwise
an explicit Runge--Kutta
method may be employed. The two subflows are then suitably
combined to produce an approximation of the exact
solution of second order in time (this is the so-called Strang splitting, also
known as split-step method). 
While such an approach is computationally attractive, it is
of low accuracy when compared, for instance, to the standard
explicit fourth-order Runge--Kutta method.
An obvious remedy is then to look for higher-order splitting methods,
or alternative schemes that perform equally well.
In theory, it is indeed possible to construct splitting schemes of arbitrary
order by suitably composing and combining the subflows,
and we refer to Reference~\cite{TCN09} for a comparison of
high-order splitting methods with real coefficients for
nonlinear Schr\"odinger equations, and to Reference~\cite{HO09b} for the
construction of high-order schemes with complex coefficients
for parabolic problems.
In our context, however, the presence of the complex term 
$\alpha_1+\rmi \beta_1$ in front of the Laplacian operator brings additional
troubles. In fact, in the general case, both high-order splitting schemes with 
real and complex coefficients
would introduce time steps with negative real parts, and thus
a potential unstable behavior since the CGL equation is irreversible.
A possible solution is the use
of additive-type splitting methods, which
employ only positive time steps (see References~\cite{D01,DLRSDLV16}).
In the very recent paper \cite{RDL24}, these schemes have been
tested (also) on the one-dimensional cubic-quintic CGL equation and proved
to be superior to classical high-order schemes.

An alternative solution is to rely on different time approximation techniques,
that still enjoy favorable properties in terms of stability.
In this sense, exponential integrators of Lawson type~\cite{L67},
also known as Integrating Factor (IF) methods,
provide a valuable alternative, since
they can be designed of arbitrary high order without incurring into
negative time
steps. These kinds of integrators require the computation of the action of
the matrix exponential and are employed in many instances in the literature,
see, for example, References~\cite{BFMMTP16,BDLV17,CEM20}.
However, to the best of our knowledge, they are of little usage in the
community of researchers interested in simulating CGL equations 
(see preprint~\cite{DK17} in which some Lawson methods are compared
with other exponential integrators and splitting schemes on one- or
two-dimensional examples with periodic boundary conditions).

In this paper, we aim at proposing efficient multidimensional
implementations 
of exponential-type schemes (i.e., of splitting and Lawson type) for
CGL equations, 
and to apply the methods to physically relevant examples such as the evolution 
of two- and three-dimensional autooscillating fields in the presence
of walls and corners,
evolution of necklace-ring patterns, and
dynamics of dissipative solitons. To accomplish the result,
we will in particular employ the
$\mu$-mode approach (see References~\cite{CCEOZ22,CCZ23}) or
the Fast Fourier Transform (FFT), depending
on the specific example under consideration, for the efficient evaluation 
of the actions of matrix exponentials needed in the linear subflow or in the
stages of the time marching methods.
{\color{black}Other numerical integrators for the CGL
  equations
which do not suffer from a time step size restriction exist in the literature.
  They involve implicit or exponential methods in time. Among them,
  we mention IMplicit EXplicit (IMEX) methods, which require
  the solution of linear systems (see the recent Reference~\cite{WLL23}),
  and exponential Runge--Kutta integrators,
  which need the approximation of matrix functions related to the exponential
  (the so-called $\varphi$-functions,
  see References~\cite{DK17,MB20}). We do not consider them in this
  manuscript since, in general, it is not possible to directly apply the
  $\mu$-mode approach to the solution of linear systems and to the
  approximation of matrix $\varphi$-functions arising from
  the discretization of multidimensional CGL equations.}

The remaining part of the manuscript is organized as follows.
In Section~\ref{sec:td} we present the schemes that will be employed for 
the time evolution of the CGL equations. We proceed in
Section~\ref{sec:expcomp} by describing in detail how to efficiently compute the
action of the linear operator and of the matrix exponential in the case of
a finite difference spatial discretization or a Fourier pseudospectral one. 
We then perform  in Section~\ref{sec:numexp} extensive numerical experiments on 
two- and three-dimensional
(coupled) CGL equations with cubic and cubic-quintic nonlinearities,
showing the 
effectiveness of the proposed exponential-type time marching methods. We finally
draw the conclusions in Section~\ref{sec:conc}.

\section{Time discretization}\label{sec:td}
As mentioned in the introduction, although explicit methods may suffer from a 
time step size restriction when employed in stiff ODEs systems of
type~\eqref{eq:ODEs},
they are usually very simple to implement and can be quite efficient,
especially if
the number of discretization points $N$, and thus the norm of the matrix
$K$, is not very large (see, for instance, References~\cite{EHGP01,RCK08}).
In our context, we consider popular Runge--Kutta {\color{black}
  methods of order two and four as
terms of comparison to the exponential-type schemes which we will propose.
The employed integrators are
\begin{equation}\label{eq:RK2}
  \begin{aligned}
    \bof_{n1}&=K\bu_n+\bg(\bu_n),\\
    \bof_{n2}&=K\left(\bu_n+\tau\bof_{n1}\right)+
    \bg\left(\bu_n+\tau\bof_{n1}\right),\\
\bu_{n+1}&=\bu_n+\frac{\tau}{2}(\bof_{n1}+\bof_{n2}),
\end{aligned}
  \end{equation}
known as Heun's method or explicit trapezoidal method, and
}
\begin{equation}\label{eq:RK4}
  \begin{aligned}
    \bof_{n1}&=K\bu_n+\bg(\bu_n),\\
    \bof_{n2}&=K\left(\bu_n+\frac{\tau}{2}\bof_{n1}\right)+
    \bg\left(\bu_n+\frac{\tau}{2}\bof_{n1}\right),\\
    \bof_{n3}&=K\left(\bu_n+\frac{\tau}{2}\bof_{n2}\right)+
    \bg\left(\bu_n+\frac{\tau}{2}\bof_{n2}\right),\\
    \bof_{n4}&=K(\bu_n+\tau\bof_{n3})+
    \bg(\bu_n+\tau\bof_{n3}),\\
\bu_{n+1}&=\bu_n+\frac{\tau}{6}\left(\bof_{n1}+2\bof_{n2}+2\bof_{n3}+\bof_{n4}\right).
\end{aligned}
  \end{equation}
{\color{black}For simplicity of presentation,
  we consider a constant time step size $\tau=t_{n+1}-t_n$,
  while $\bu_{n}$ and  $\bu_{n+1}$ denote
  the approximations of $\bu(t_n)$ and $\bu(t_{n+1})$, respectively}.

Concerning the exponential-type methods, we will consider two families:
the splitting and the Lawson ones.

\subsection{Splitting schemes}
The most popular time integrator for systems of ODEs~\eqref{eq:ODEs} arising
from the spatial discretization of CGL equations
is probably the Strang splitting scheme
(see, for instance, 
References~\cite{HFDW06,ZBD07,ZY15b,DK18}).
If we denote by
\begin{equation*}
\check\bu(t)=\phi_t(\check\bu_0)
\end{equation*}
the flow of $\check\bu'(t)=\bg(\check\bu(t))$ starting from
$\check \bu_0$,
then the Strang splitting scheme {\color{black}for system~\eqref{eq:ODEs}}
is defined by
\begin{equation}\label{eq:Strang}
\bu_{n+1}=\mathcal{S}_\tau(\bu_n)=
  \phi_{\tau/2}\left(\rme^{\tau K}\phi_{\tau/2}(\bu_n)\right).
\end{equation}
It is a method of second order with respect to the time step size $\tau$.
The linear flow is the action of the matrix exponential $\rme^{\tau K}$.
The nonlinear flow $\phi_{\tau/2}$ can be computed exactly
(e.g., for the cubic or the quintic CGL equations, see
Section~\ref{sec:numexp}) or conveniently approximated by
{\color{black} explicit methods
(the Runge--Kutta schemes presented above, for instance)},
since the stiffness is already handled by the linear flow.
A related method is presented in Reference~\cite{WZ13}, where
the linear flow is treated in a different way. In fact,
it is first observed that the flow corresponding to
\begin{equation*}
  \partial_t \hat u=(\alpha_1+\rmi\beta_1)\Delta \hat u
\end{equation*}
can be \emph{exactly} split into the partial flows of
\begin{equation*}
  \partial_t \hat u_\mu=
  (\alpha_1+\rmi\beta_1)\partial_{x_\mu x_\mu}\hat u_\mu,\quad \mu=1,\ldots,d,
\end{equation*}
where $\partial_{x_\mu x_\mu}$ is discretized by finite differences.
Then, each of the discrete partial linear flows
(i.e., the action of the matrix
exponential)
is approximated by the Crank--Nicolson
scheme.
{\color{black}Similarly,
  in Reference~\cite{W10} the Crank--Nicolson method is applied
  to a Chebyshev tau spectral discretization in space.}
Instead of the second-order approximation in time
of the discretized flow given by the underlying trapezoidal rule,
we will see in Section~\ref{sec:expcomp} how it is possible to
evaluate \emph{efficiently and exactly}
the action of the discretized flow $\rme^{\tau K}$
in equation~\eqref{eq:Strang}.

Standard
composition methods involving the flows $\phi_{a_j\tau}$ and $\rme^{b_j\tau K}$
can generate higher-order methods. However, as mentioned previously in the
introduction, they would require some of
the coefficients $a_j$ and $b_j$ to have real negative part or complex
conjugate imaginary parts with large arguments,
making them ill-suited for irreversible problems
like the CGL equations
(see, for instance,
Reference~\cite{CCDV09,HO09b,DLRSDLV16}). Therefore, for our purposes
we will consider the following fourth-order splitting scheme
\begin{equation}\label{eq:richextrap}
\bu_{n+1} = \frac{4}{3}\mathcal{S}_{\tau/2}(\mathcal{S}_{\tau/2}(\bu_n))-
  \frac{1}{3}\mathcal{S}_{\tau}(\bu_n),
\end{equation}
which can be explicitly written as
\begin{equation}\label{eq:Richardson}
  \begin{aligned}
    \bu_{n2}&=\phi_{\tau/2}\left(\rme^{\tau K}\phi_{\tau/2}(\bu_n)\right),\\
    \bu_{n+1}&=\frac{4}{3}\phi_{\tau/4}\left(\rme^{\tau K/2}\phi_{\tau/2}\left(\rme^{\tau K/2}\phi_{\tau/4}(\bu_n)\right)\right)
    -\frac{1}{3}\bu_{n2}.
  \end{aligned}
\end{equation}
This scheme can be seen as the Richardson extrapolation of the Strang splitting 
method, and was developed and analyzed in Reference~\cite{D01}.
It is an additive scheme similar to the so-called affine splitting schemes
introduced in Reference~\cite{DLRSDLV16}, but compared to those methods it requires
less flow evaluations. Another splitting scheme of order four
with positive steps in the linear flow $\rme^{\tau K}$ is presented
in the conclusions of Reference~\cite{CCDV09}, but in fact it requires
nine flow evaluations. For this reason, in this work we limit ourselves to the
scheme~\eqref{eq:richextrap} among high-order splitting methods.

\subsection{Lawson schemes}
Lawson methods
were presented in Reference~\cite{L67} as generalized
Runge--Kutta methods for systems with large Lipschitz constants, hence
well-suited for many stiff problems.
They can be derived by
{\color{black}formally rewriting  system~\eqref{eq:ODEs} in terms of the
new unknown $\bw(t) = \rme^{-t K}\bu(t)$, where $\rme^{-t K}$
  acts as an integrating factor.  Applying an explicit
  Runge--Kutta method to the new system
  $\bw'(t)=\rme^{-tK}\bg (\rme^{tK}\bw(t))$,
  and rewriting the scheme in terms of the original
  unknown, the desired integrator is obtained.}
For instance, if we apply Heun's method \eqref{eq:RK2}
to the differential equation for $\bw(t)$, we
obtain the Lawson method of order two
known in the literature as Lawson2b
\begin{equation}\label{eq:lawson2b}
\begin{aligned}
  \bu_{n2}&=\rme^{\tau K}(\bu_n+\tau \bg(\bu_n)),\\
  \bu_{n+1}&=\rme^{\tau K}\left(\bu_n+\frac{\tau}{2} \bg(\bu_n)\right)+
  \frac{\tau}{2}\bg(\bu_{n2}).
\end{aligned}
\end{equation}
If we apply instead the Runge--Kutta method of order
four~\eqref{eq:RK4} to the system of ODEs for $\bw(t)$,
we obtain the Lawson scheme of order four
\begin{equation}\label{eq:lawson4}
\begin{aligned}
  \bu_{n2}&=\rme^{\tau K/2}\left(\bu_n+\frac{\tau}{2} \bg(\bu_n)\right),\\
  \bu_{n3}&=\rme^{\tau K/2}\bu_n+\frac{\tau}{2}\bg(\bu_{n2}),\\
  \bu_{n4}&=\rme^{\tau K}\bu_n+\tau\rme^{\tau K/2}\bg(\bu_{n3}),\\
  \bu_{n+1}&=\rme^{\tau K}\left(\bu_n+\frac{\tau}{6} \bg(\bu_n)\right)+
  \frac{\tau}{3}\rme^{\tau K/2}(\bg(\bu_{n2})+\bg(\bu_{n3}))+
  \frac{\tau}{6}\bg(\bu_{n4}),
\end{aligned}
\end{equation}
which is typically referred to as Lawson4.
Similarly to the splitting schemes, they require only the evaluation
of the action of the matrix exponential, and in particular not of the 
$\varphi$-functions needed by other
classes of exponential methods \cite{HO10}.
Notice finally that, by construction, these methods do not have
negative
coefficients in front of the time step size in the matrix exponentials.

\section{Action of the linear
  operator and of its matrix exponential}\label{sec:expcomp}
\subsection{Finite difference spatial discretization}\label{sec:fd}
When system~\eqref{eq:ODEs} is the result of a finite difference approximation
in space {\color{black}with homogeneous Dirichlet or Neumann boundary conditions,
the matrices $A_\mu$ are sparse with a certain number
of non-zero diagonals adjacent to the main one (the specific number depends on
the chosen approximation).}
The time
marching schemes described in the previous section require, in particular,
the action of the operator $K$ and of its (scaled) matrix exponential, whose 
efficient computation is crucial for the effectiveness of the overall approach.
To this aim, it is possible to suitably exploit the Kronecker sum
structure of the linear operator to provide equivalent tensor formulations
of the tasks. A reader not familiar with the following
tensor formalism is invited to 
check Reference~\cite{CCZ23} for a detailed explanation.
In fact, we can compute the action of the matrix 
$K\in\CC^{N\times N}$ on a vector $\bu\in\CC^{N}$ by exploiting
the equivalence
\begin{equation}\label{eq:kronsumv}
K\bu=\mathrm{vec}\left(\sum_{\mu=1}^d \bU\times_\mu A_\mu\right),
\end{equation}
without assembling the matrix $K$ itself. Here,
$\bU\in\CC^{n_1\times\cdots\times n_d}$
is the order-$d$ tensor such that $\mathrm{vec}(\bU)=\bu$,
where $\mathrm{vec}$ is the operator which stacks the columns of the
input tensor into a suitable column vector. Moreover, $\times_\mu$
denotes the so-called $\mu$-mode product. This operation
multiplies the matrix $A_\mu$
onto the $\mu$-fibers (i.e., the generalizations
of rows and columns to tensors) of the tensor $\bU$. The action of the
matrix exponential $\rme^{\tau K}$ can also be computed using $\mu$-mode
products, again without assembling the matrix $K$. Indeed, it holds
\begin{equation}\label{eq:Tucker}
  \rme^{\tau K}\bu=\mathrm{vec}\left(\bU\times_1\rme^{\tau A_1}\times_2
  \cdots\times_d \rme^{\tau A_d}\right).
\end{equation}
The concatenation of $\mu$-mode products is usually referred to
as a
Tucker operator, and costs $\mathcal{O}((n_1+\cdots+n_d)N)$
{\color{black}floating-point} operations.
Moreover, after computing the relevant matrix exponentials,
it can be performed using the high performance
level 3 BLAS, since the matrices $\rme^{\tau A_\mu}$ and the tensor $\bU$
are \emph{dense}. This approach leads to a huge computational saving,
see Reference~\cite{CCEOZ22}.
{\color{black}Notice that formula~\eqref{eq:Tucker} is
  peculiar to the exponential and does not hold, for instance, for
  the inverse function (needed by implicit methods)
  or the $\varphi$-functions (needed by other exponential integrators).
  Also, we remark that each matrix exponential of small
  size $n_\mu\times n_\mu$ can be computed
by using standard Pad\'e (see Reference~\cite{AMH10})
or Taylor approximations (see References~\cite{CZ19,SID19,AIDAJ23}),
with a cost $\mathcal{O}(n_\mu^3)$.
These computations can be
performed once and for all before the time integration starts if
the time step size is constant. Otherwise,
the matrix exponentials must be
recomputed at each time step.}
We finally highlight that the sparsity pattern of the
matrices $A_\mu$ is \emph{not} important. In fact, they are of small size, and
the resulting matrix exponentials are full. Hence, there is no considerable 
advantage in performing computations with sparse matrices. 
The limit case would be, for instance, the usage of dense spectral 
differentiation matrices, see in particular References~\cite{B00,CCEO24}.
Also, the substitution of the Laplace operator in equation~\eqref{eq:cqCGL} 
with its fractional counterpart 
(i.e., considering a space-fractional CGL equation)
would lead to full discretization matrices (see, for example,
References~\cite{ZOG21,ZZS23}).

\subsection{Pseudospectral Fourier spatial discretization}\label{sec:Fourier}
When dealing with periodic boundary conditions,
it is common to employ a pseudospectral Fourier
decomposition in space. We consider without loss of
generality the spatial domain {\color{black}$\Omega=(0,2\pi)^d$}.
The generalization to arbitrary
domains {\color{black}$(a_1,b_1)\times\cdots\times(a_d,b_d)$}
is straightforward. Then,
we can approximate the solution of the CGL equation~\eqref{eq:cqCGL}
by
\begin{equation*}
  u(t,\bx)\approx \sum_{j_d=1}^{n_d}\cdots \sum_{j_1=1}^{n_1}
  u_{j_1\ldots j_d}(t)\phi_{j_1}(x_1)\cdots\phi_{j_d}(x_d),
\end{equation*}
where
\begin{equation*}
\phi_{j_\mu}(x_\mu)=\rme^{\rmi (j_\mu-1) x_\mu}
\end{equation*}
and $\bu(t)=(u_{j_1\ldots j_d}(t))$ are the discrete
Fourier coefficients of $u(t,\bx)$.
By performing standard calculations, we then end up with a system of ODEs for
the spectral coefficients in the form~\eqref{eq:ODEs}, where
the matrices $A_\mu$
are diagonal and defined by
\begin{equation*}
  A_\mu=\mathrm{diag}\left(
  \frac{\alpha_2}{d},
  (\alpha_1+\rmi\beta_1)(-1)+  \frac{\alpha_2}{d},\ldots,
  (\alpha_1+\rmi\beta_1)(-(n_\mu-1)^2)+  \frac{\alpha_2}{d}\right).
\end{equation*}
The nonlinear function $\bg(\bu(t))$ is in this case
\begin{equation*}
\mathcal{F}(g(\mathcal{F}^{-1}(\bu(t)))),
\end{equation*}
where $\mathcal{F}$ denotes the $d$-dimensional
Fourier transform (which maps a function
sampled at a Cartesian grid of size $n_1\times\cdots\times n_d$
in {\color{black}$(0,2\pi)^d$} 
to its discrete Fourier coefficients). {\color{black}It can be realized by
  the Fast Fourier Transform (FFT) algorithm
  at a cost $\mathcal{O}((\log n_1+\cdots+\log n_d)N)$.}
In this setting, the action of $K$ on a vector $\bu$ can be easily computed as
the pointwise multiplication between the tensors
\begin{equation*}
  \boldsymbol K=
  (k_{j_1\ldots j_d})=
    \left((\alpha_1+\rmi\beta_1)\left(-\sum_{\mu=1}^d(j_\mu-1)^2\right)+\alpha_2\right)
\end{equation*}
and $\bU$, where $\mathrm{vec}(\bU)=\bu$. Similarly, since
$\rme^{\tau A_\mu}$ are diagonal matrices, the action of the
matrix exponential $\rme^{\tau K}$ on the vector $\bu$
can be computed as the pointwise multiplication between the tensors
\begin{equation*}
  \boldsymbol E=
  (e_{j_1\ldots j_d})=
    \left(\rme^{\tau ((\alpha_1+\rmi\beta_1)(-\sum_{\mu=1}^d(j_\mu-1)^2)+\alpha_2)}\right)
\end{equation*}
and $\bU$. {\color{black}Forming the tensor $\boldsymbol E$
costs $\mathcal{O}(N)$ and, similarly to the finite difference case, can
  be done once and for all before the time integration starts if
  the time step size $\tau$ is constant.}

We finally remark here that other pseudospectral decompositions could profit
from the underlying tensor structure of the problem. For instance, 
another classical pseudospectral decomposition employed for Schr\"odinger or
CGL equations is based on Hermite functions
(see References~\cite{CCEOZ22,RDL24}), which naturally vanish for
$\bx\to\infty$. In this case, the resulting system of ODEs for the spectral
coefficients is not diagonal, but still possesses Kronecker sum structure 
(the discretization matrices $A_\mu$ are full with nonzero even diagonals).
Hence, we can still efficiently compute the action of the linear operator $K$
and its matrix exponential by means of the $\mu$-mode approach
described in Section~\ref{sec:fd}.

\section{Numerical experiments}\label{sec:numexp}
In this section, we perform some numerical experiments with a variety of
CGL equations {\color{black}defined on Cartesian product domains}
to assess the efficiency of high-order
exponential-type methods over other popular techniques.

{\color{black}
In all the examples, we first compare the performances of the methods in
reaching a range of accuracies (measured in the infinity norm
relatively with respect to a reference solution).
To this aim, we choose different sequences of time steps and measure the overall
needed wall-clock time of the integration. It turns out that, for each
numerical experiment,
 the same number of time steps can be employed
 for second-order and fourth-order exponential-type methods,
 respectively. In addition, the selected
 sequences allow to observe the expected orders of convergence.
 Concerning the Runge--Kutta methods, we choose the number of time steps
 also to clearly show the need of a time step size restriction.}

Then, in a second experiment we show the dynamics of the equation up
to a large final time. In this case, the time marching is performed with
a single time step size employing the integrator that performed best in the
previous experiment.

All the numerical examples are done on a standard laptop equipped with an
Intel\textsuperscript{\textregistered} Core\textsuperscript{\texttrademark}
i7-10750H chip (six physical cores) and 16GB of RAM. As programming language
we employ MATLAB (interpreted by the software
MathWorks MATLAB\textsuperscript{\textregistered} R2022a).
The source code to reproduce all the experiments, fully compatible with
GNU Octave, can be found in a GitHub
repository\footnote{Available at \url{https://github.com/caliarim/expCGL}.}.

In all cases, we employ the tensor formulation of the problem using the 
built-in MATLAB N-D array type as storage for the tensors (i.e.,
no external class is needed).
In particular, when performing experiments with finite differences, both the
actions of the linear operator $K$ and
of the matrix exponentials are performed, in a $\mu$-mode fashion,
by using the right hand sides of
formulas~\eqref{eq:kronsumv} and~\eqref{eq:Tucker}, without vectorizing the
results.
In practice, the tensor-matrix operations are done 
with level 3 BLAS using the functions \texttt{kronsumv} and \texttt{tucker}
of KronPACK\footnote{Available at \url{https://github.com/caliarim/KronPACK}.}
(see Reference~\cite{CCZ23}).
The needed small sized matrix exponentials are
computed using a rational Pad\'e approach with scaling and
squaring~\cite{AMH10} (implemented in the MATLAB function \texttt{expm}).
When performing experiments with periodic boundary conditions, we
exploit the high performance FFTW library
(through the MATLAB function \texttt{fftn})
to realize a Fourier pseudospectral
discretization, and all the relevant operations are performed pointwise
(see Section~\ref{sec:Fourier}).
{\color{black}We work in a constant time step size scenario.
  Hence, the computation of the matrix exponentials is
  performed once and for all before the time integration starts,
as mentioned in the previous section.
This phase has negligible computational cost compared to that of the time
integration, see also
the discussions in References~\cite{CCEOZ22,CC24,C24} for more
details. The investigation of reliable and efficient variable time step
size methods is beyond the scope of the present manuscript.}

Finally, for sake of compact notation, we set some labels for the time marching
methods introduced in Section~\ref{sec:td}. More in detail, the Strang splitting
scheme~\eqref{eq:Strang} is denoted as \textsc{strang},
the Lawson2b integrator~\eqref{eq:lawson2b} as \textsc{if2},
{\color{black}the second-order Runge--Kutta method~\eqref{eq:RK2} as
\textsc{rk2},}
the fourth-order splitting method~\eqref{eq:Richardson} 
as \textsc{split4},
the Lawson4 scheme~\eqref{eq:lawson4} as \textsc{if4},
and the standard explicit fourth-order Runge--Kutta integrator~\eqref{eq:RK4} 
as \textsc{rk4}.
\subsection{2D cubic CGL equation}\label{sec:CCGL2D}
In the first numerical example, we consider the cubic CGL equation
\begin{equation}\label{eq:ccgl}
  \partial_tu=(\alpha_1+\rmi\beta_1)\Delta u+
  \alpha_2 u+
  (\alpha_3+\rmi\beta_3)\lvert u\rvert^2u,
\end{equation}
in the two-dimensional domain {\color{black}$\Omega=(0,100)^2$}
with the parameters $\alpha_1=1$, $\beta_1=2$, $\alpha_2=1$, $\alpha_3=-1$,
and $\beta_3=0.2$. The number of discretization points is equal to 
$n=256$ in both directions (total number of degrees of freedom
$N=256^2$).
The values of the initial datum are random of \textit{small amplitude}, that is
taken from a standard normal distribution divided by 5000. In practice, this is
realized in MATLAB by the command \verb+randn(n,n)/5000+ with the random number
generator specified by \verb+rng('default')+.
For this example, we will consider two different types of boundary conditions 
(i.e., homogeneous Dirichlet and periodic), which lead to substantially 
different dynamics.
The parameters, the initial datum, and the boundary conditions are taken from
Reference~\cite{EHGP01}.

Concerning the time integration methods of splitting type (i.e., \textsc{strang}
and \textsc{split4}), we notice that the nonlinear part has exact solution.
Indeed, the analytical solution of
\begin{equation*}
  \check{\bu}'(t) = (\alpha_3+\rmi\beta_3)\lvert \check{\bu}(t)\rvert^2\check{\bu}(t)
\end{equation*}
starting from the initial datum $\check{\bu}_0$ is
\begin{equation}\label{eq:exsolcub}
  \check{\bu}(t) = \phi_t^{\mathrm c}(\check{\bu}_0)=\rme^{-\frac{\alpha_3+\rmi\beta_3}{2\alpha_3}
                    \log\lvert1-2\alpha_3\lvert \check{\bu}_0 \rvert^2t\rvert}
                    \check{\bu}_0.
\end{equation}

\subsubsection{Homogeneous Dirichlet boundary conditions}\label{sec:CCGL2Dhd}

We first consider the example with homogeneous Dirichlet boundary
conditions.
We then approximate the Laplace operator by fourth-order 
finite differences using internal points. The resulting unidirectional
discretization matrix is
\begin{equation*}
  \frac{1}{12h^2}\begin{bmatrix}
    -15 & -4  & 14 & -6 & 1 & 0 & \ldots & \ldots & 0\\
    16 & -30 & 16 & -1 & 0 & \ldots & \ldots & \ldots & 0\\
    -1 & 16 & -30 & 16 & -1 & 0 & \ldots & \ldots & 0\\
    0 & -1 & 16 & -30 & 16 & -1 & 0 & \ldots & 0\\
    \vdots & \ddots & \ddots & \ddots & \ddots & \ddots &\ddots & \ddots & \vdots\\
    0 & \ldots & 0 & -1 & 16 & -30 & 16 & -1 & 0\\
    0 & \ldots & \ldots & 0 & -1 & 16 & -30 & 16 & -1\\
    0 & \ldots & \ldots & \ldots & 0 & -1 & 16 & -30 & 16\\
    0 & \ldots & \ldots & 0 & 1 & -6 & 14 & -4 & -15
  \end{bmatrix}
\end{equation*}
for both directions, where $h=100/(n+1)$ denotes the spatial step size. 
Notice that in the first and in the last
rows (which correspond to the first and last internal point of 
the discretized domain, respectively) we employed non-centered 
finite differences to overcome the lack of enough surrounding points for the
usual centered stencil.

\begin{table}[!htb]
  \centering
  {\small
  \begin{tabular}{llll|llll}
    \hline
    \multicolumn{4}{l|}{\textsc{if2}} &                    \multicolumn{4}{l}{\textsc{if4}} \\
    steps & time (s) & error & order &                 steps & time (s) & error & order \\
    \hline                                           
    25 & $2.245\rme{-01}$ & $4.359\rme{-07}$ &  ---  &  100 & $1.669\rme{+00}$ & $5.858\rme{-10}$ & ---   \\
    75 & $3.386\rme{-01}$ & $4.984\rme{-08}$ &  1.974 & 150 & $2.543\rme{+00}$ & $2.378\rme{-10}$ & 2.223 \\
    125 & $5.389\rme{-01}$ & $1.842\rme{-08}$ & 1.948 & 200 & $5.884\rme{+00}$ & $1.033\rme{-10}$ & 2.899 \\
    175 & $7.183\rme{-01}$ & $9.538\rme{-09}$ & 1.957 & 250 & $7.579\rme{+00}$ & $4.800\rme{-11}$ & 3.433 \\
    225 & $9.734\rme{-01}$ & $5.782\rme{-09}$ & 1.992 & 300 & $1.132\rme{+01}$ & $2.375\rme{-11}$ & 3.860 \\
    \hline
    \hline
    \multicolumn{4}{l|}{\textsc{strang}}&                                 \multicolumn{4}{l}{\textsc{split4}}\\
    steps & time (s) & error & order&                           steps & time (s) & error & order\\
    \hline                                                                                                                   
    25 & $1.161\rme{-01}$ & $4.359\rme{-07}$ & ---     & 100 & $9.419\rme{-01}$& $5.858\rme{-10}$ & ---     \\
     75 & $2.362\rme{-01}$ & $4.984\rme{-08}$ & 1.974   & 150 & $1.399\rme{+00}$ & $2.378\rme{-10}$ & 2.223 \\
     125 & $3.015\rme{-01}$ & $1.842\rme{-08}$ & 1.948  & 200 & $3.036\rme{+00}$ & $1.033\rme{-10}$ & 2.899 \\
     175 & $3.575\rme{-01}$ & $9.538\rme{-09}$ & 1.957  & 250 & $3.752\rme{+00}$ & $4.797\rme{-11}$ & 3.437 \\
     225 & $5.012\rme{-01}$ & $5.782\rme{-09}$ & 1.992  & 300 & $5.429\rme{+00}$ & $2.368\rme{-11}$ & 3.872 \\
    \hline
    \hline
    \multicolumn{4}{l|}{\textsc{rk2}} &     \multicolumn{4}{l}{\textsc{rk4}}\\
        steps & time (s) & error & order & steps & time (s) & error & order \\
    \hline
    355 & $\times$ & $\times$ & --- &
    255 & $\times$ & $\times$ & --- \\
    455 & $\times$ & $\times$ & $\times$ &
    355 & $3.382\rme{+00}$ & $3.241\rme{-04}$ & $\times$\\
    555 & $2.179\rme{+00}$ & $1.088\rme{-04}$ & $\times$  &
    455 & $4.086\rme{+00}$ & $1.467\rme{-09}$ & 49.58\\
    655 & $2.599\rme{+00}$ & $7.821\rme{-05}$ & 1.993 &
    555 & $4.608\rme{+00}$ & $6.594\rme{-10}$ & 4.024\\
    755 & $2.956\rme{+00}$ & $5.891\rme{-05}$  & 1.994  &
    655 & $5.587\rme{+00}$ & $3.396\rme{-10}$ & 4.006\\
    \hline
  \end{tabular}%
  }
  \caption{Number of time steps, wall-clock time (in seconds),
      relative error at final time $T=6$, and observed numerical order of convergence
      for the solution of the 2D cubic CGL equation~\eqref{eq:ccgl} with homogeneous
      Dirichlet boundary conditions and different integrators. 
      The total number of degrees of freedom is $N=256^2$.
      In the table, the cross symbol $\times$ means NaN/Inf value returned.
      See also Figure~\ref{fig:cgl_cubic_2D_FD_cpudiag} for a graphical representation.}
  \label{tab:cgl_cubic_2D_FD}
\end{table}
\begin{figure}[!htb]
  \centering
%
%
%
\begin{tikzpicture}

\begin{axis}[%
width=2.6in,
height=1.7in,
at={(0.758in,0.481in)},
scale only axis,
xmode=log,
xmin=0.06,
xmax=15,
xminorticks=true,
xlabel style={font=\color{white!15!black}},
ymode=log,
ymin=1e-11,
ymax=1e-03,
yminorticks=true,
ylabel style={font=\color{white!15!black}},
xlabel = {Wall-clock time (s)},
ylabel = {Relative error},
axis background/.style={fill=white},
legend style={at={(0.025,0.65)}, anchor=south west, legend cell align=left, align=left, draw=white!15!black, font=\footnotesize},
legend columns = 2
]
\addplot [color=blue, line width=1pt, mark size = 2pt, mark=triangle, mark options={solid, rotate=180, blue}]
  table[row sep=crcr]{%
0.224468	4.35903461509076e-07\\
0.338607	4.98376006841711e-08\\
0.53893	1.84240075363468e-08\\
0.718266	9.53795389022339e-09\\
0.973402	5.7821066343833e-09\\
};
\addlegendentry{\textsc{if2}}

\addplot [color=cyan, line width=1pt, mark size = 2pt,mark=square, mark options={solid, cyan}]
  table[row sep=crcr]{%
1.66906	5.85820173942779e-10\\
2.542625	2.3781230154309e-10\\
5.8841	1.03275899585481e-10\\
7.579321	4.80045894532108e-11\\
11.31589	2.3747540908261e-11\\
};
\addlegendentry{\textsc{if4}}

\addplot [color=red, line width=1pt, mark size = 2pt,mark=o, mark options={solid, red}]
  table[row sep=crcr]{%
0.116127	4.35893610044074e-07\\
0.236159	4.98372771530963e-08\\
0.301468	1.84239410204334e-08\\
0.357455	9.5379324532401e-09\\
0.501239	5.78209796791114e-09\\
};
\addlegendentry{\textsc{strang}}

\addplot [color=magenta, line width=1pt, mark size = 2pt,mark=asterisk, mark options={solid, magenta}]
  table[row sep=crcr]{%
0.941883	5.85811471005438e-10\\
1.399451	2.37812418099952e-10\\
3.036496	1.03298290388486e-10\\
3.752489	4.79691184677429e-11\\
5.429336	2.36814570825637e-11\\
};
\addlegendentry{\textsc{split4}}

\addplot [color=black, line width=1pt, mark size = 2pt,mark=triangle, mark options={solid, black}]
  table[row sep=crcr]{%
2.1788e+00	1.0881e-04\\
2.5993e+00	7.8211e-05\\
2.9563e+00	5.8913e-05\\
};
\addlegendentry{\textsc{rk2}}

\addplot [color=green, line width=1pt, mark size = 2pt,mark=diamond, mark options={solid, green}]
  table[row sep=crcr]{%
3.381788	0.000324087588321909\\
4.086068	1.46668893292526e-09\\
4.607702	6.59424979316862e-10\\
5.586881	3.39603561263342e-10\\
};
\addlegendentry{\textsc{rk4}}

\end{axis}

\end{tikzpicture}%
  \caption{Results for the simulation of the 2D cubic CGL
    equation~\eqref{eq:ccgl} with homogeneous Dirichlet
    boundary conditions and
  $N=256^2$ spatial discretization points. The number
  of time steps for each integrator is reported in
  Table~\ref{tab:cgl_cubic_2D_FD}. The final simulation time is $T=6$.}
  \label{fig:cgl_cubic_2D_FD_cpudiag}
\end{figure}
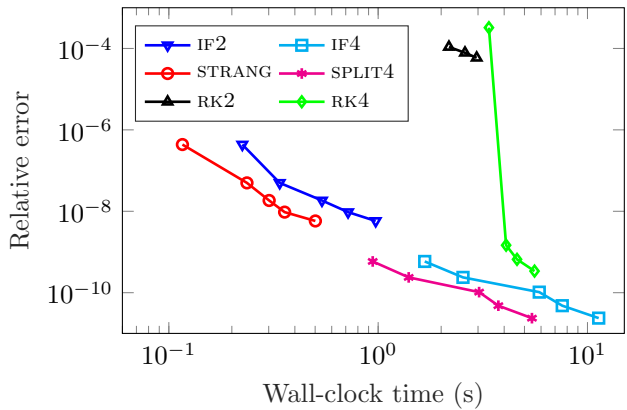

For the first experiment we set the final simulation time to $T=6$.
The detailed outcome, together with the numbers of time steps and
the obtained numerical order of convergence for the different integrators,
is reported in Table~\ref{tab:cgl_cubic_2D_FD}.
The results are also graphically depicted in a precision-work diagram in
Figure~\ref{fig:cgl_cubic_2D_FD_cpudiag}.
First of all, for sufficiently small time step sizes all the
integrators show the expected order of convergence.
{\color{black} We notice that the second-order
  exponential-type methods seem to achieve
  the same accuracy. This is true also for the fourth-order
  exponential-type schemes. However,
  the methods do not exactly compute the same
  numerical solution. Indeed,
  if we perform a simulation with $T=30$ and 50 time steps,
  the relative difference in infinity norm between
  the second-order methods and the fourth-order ones
  is $5.834\rme{-01}$ and $1.372\rme{-01}$,
  respectively.}
The explicit Runge--Kutta {\color{black}methods \textsc{rk2} and \textsc{rk4} are
the only schemes which suffer} from a restriction on the 
{\color{black}time step size}. Indeed,
as expected a stability constraint is required due to the underlying stiffness of the
problem. In fact, {\color{black} \textsc{rk2} is the least efficient scheme.
The \textsc{rk4} integrator needs more than five times the number of time steps
employed by the fourth-order exponential-type methods to reach comparable 
accuracies, resulting then less efficient in the precision-work diagram.}
Notice also that, for each integrator, the average cost per time step is 
essentially constant. This is again expected, since all the methods under consideration
are direct (see also the discussion in Reference~\cite{CC24}).
{\color{black}
Splitting methods (both of second and of fourth order) 
appear to be more efficient than the corresponding Lawson ones.
Moreover, if we extend the lines
in Figure~\ref{fig:cgl_cubic_2D_FD_cpudiag} to cover a wide range of
accuracies, we observe that the fourth-order methods are slightly
better than the second-order schemes. Overall,
\textsc{split4} allows to reach stringent accuracies
with the smallest amount of computational time.}

\begin{figure}[!htb]
  \centering
  \input{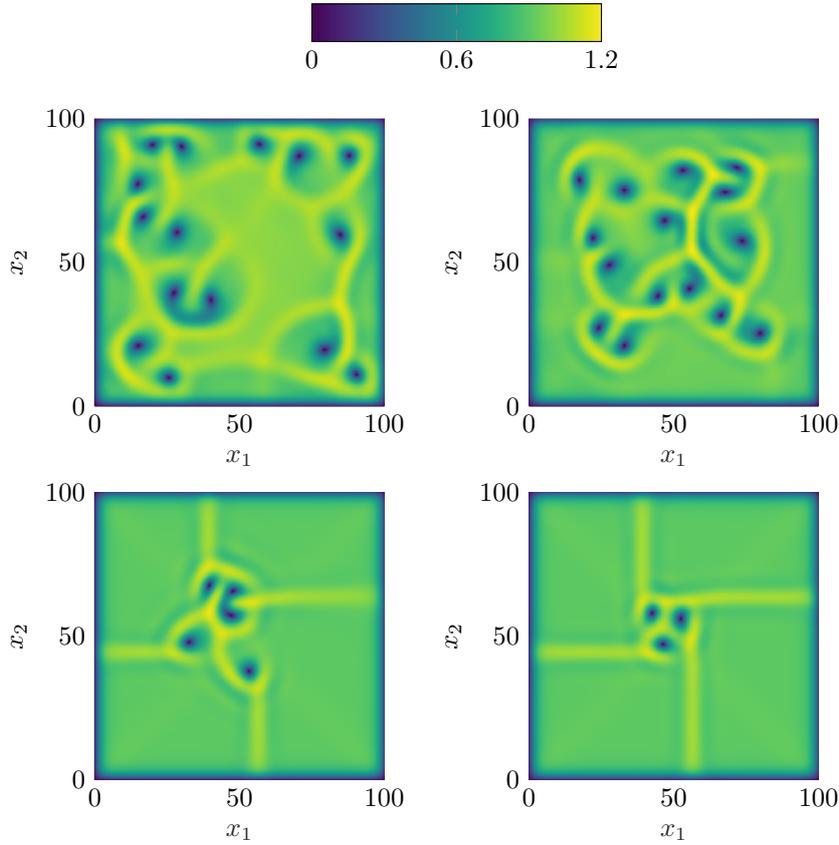}
  \caption{Dynamics of $\lvert u \rvert$ for the 2D cubic CGL
    equation~\eqref{eq:ccgl} with homogeneous Dirichlet
    boundary conditions and
  $N=256^2$ spatial discretization points. The time integrator is 
  \textsc{split4}. The final times are $T=30$ (top left), $T=60$ (top right),
  $T=120$ (bottom left), and $T=240$ (bottom right) with a number of time steps
  equal to
  $m=50$, $m=100$, $m=200$, and $m=400$, respectively. The
  wall-clock times of the simulations are $0.54$, $1.04$, $1.76$, and $3.62$ seconds, respectively.}
  \label{fig:cgl_cubic_2D_FD_dyn}
\end{figure}

Finally, in Figure~\ref{fig:cgl_cubic_2D_FD_dyn} we show the evolution of the
dynamics of the cubic CGL equation under consideration. To this aim, we 
simulate with \textsc{split4} up to final times $T=30$, $T=60$, $T=120$, 
and $T=240$ (with $m=50$, $m=100$, $m=200$, and $m=400$ time steps,
respectively).
The latter
simulation took less than 4 seconds. The results
match satisfactorily with the theoretical predictions and with 
what already reported in the literature (see in particular Reference~\cite{EHGP01}).
Indeed, starting from a random initial condition, defects are formed and they
move around the domain while being pushed to the center from the waves emitted
by the boundaries. Finally, after $T=240$, the system reaches a so-called
\textit{frozen state}, and the modulus of the solution does not change anymore.

\subsubsection{Periodic boundary conditions}
We now consider the cubic CGL equation~\eqref{eq:ccgl}
completed with periodic boundary conditions. In this case, the space discretization
is performed with a Fourier pseudospectral method.

\begin{table}[!htb]
  \centering
  {\small
  \begin{tabular}{llll|llll}
    \hline
    \multicolumn{4}{l|}{\textsc{if2}} &                      \multicolumn{4}{l}{\textsc{if4}} \\
        steps & time (s) & error & order &                  steps & time (s) & error & order \\
    \hline                                              
    25 & $6.806\rme{-02}$ & $4.322\rme{-07}$ & ---  &   200 & $1.038\rme{+00}$ & $1.846\rme{-10}$ & ---  \\
    75 & $1.742\rme{-01}$ & $5.055\rme{-08}$ & 1.953 &  250 & $1.244\rme{+00}$ & $1.012\rme{-10}$ & 2.694\\
    125 & $2.811\rme{-01}$ & $1.877\rme{-08}$ & 1.940 & 300 & $1.481\rme{+00}$ & $5.742\rme{-11}$ & 3.110\\
    175 & $3.881\rme{-01}$ & $9.809\rme{-09}$ & 1.928 & 350 & $1.807\rme{+00}$ & $3.378\rme{-11}$ & 3.441\\
    225 & $5.305\rme{-01}$ & $6.014\rme{-09}$ & 1.946 & 400 & $1.949\rme{+00}$ & $2.060\rme{-11}$ & 3.706\\
    \hline
    \hline
      \multicolumn{4}{l|}{\textsc{strang}}&                            \multicolumn{4}{l}{\textsc{split4}}\\                      
      steps & time (s) & error & order&                      steps & time (s) & error & order\\                   
      \hline                                                 
    25 & $4.930\rme{-02}$ & $4.322\rme{-07}$ & ---     &    200 & $8.717\rme{-01}$& $1.846\rme{-10}$ & ---   \\
    75 & $9.104\rme{-02}$ & $5.055\rme{-08}$ & 1.953   &      250 & $9.997\rme{-01}$ & $1.012\rme{-10}$ & 2.695 \\  
     125 & $1.461\rme{-01}$ & $1.877\rme{-08}$ & 1.940 &    300 & $1.231\rme{+00}$ & $5.742\rme{-10}$ & 3.108 \\  
     175 & $2.027\rme{-01}$ & $9.809\rme{-09}$ & 1.928 &    350 & $1.425\rme{+00}$ & $3.378\rme{-11}$ & 3.440 \\ 
     225 & $2.643\rme{-01}$ & $6.014\rme{-09}$ & 1.946 &    400 & $1.588\rme{+00}$ & $2.060\rme{-11}$ & 3.705 \\  
    \hline
    \hline
    \multicolumn{4}{l|}{\textsc{rk2}} &    \multicolumn{4}{l}{\textsc{rk4}}\\
        steps & time (s) & error & order &        steps & time (s) & error & order \\
    \hline
    755 & $\times$ & $\times$ & ---  &
    555 & $\times$ & $\times$ & --- \\
    855 & $\times$  & $\times$  & $\times$ &
    655 & $2.745\rme{+00}$ & $6.015\rme{-08}$ & $\times$\\
    955 & $1.976\rme{+00}$  & $3.732\rme{-05}$ & $\times$ &
    755 & $3.104\rme{+00}$ & $2.205\rme{-10}$ & 39.47\\
    1055 & $2.217\rme{+00}$ & $3.060\rme{-05}$  & 1.996 &
    855 & $3.937\rme{+00}$ & $1.327\rme{-10}$ & 4.089\\
    1155 & $2.434\rme{+00}$ & $2.554\rme{-05}$ & 1.996 &
    955 & $4.329\rme{+00}$ & $8.491\rme{-11}$ & 4.034\\
    \hline
  \end{tabular}%
  }
  \caption{Number of time steps, wall-clock time (in seconds),
  relative error at final time $T=6$, and observed numerical order of convergence
      for the solution of the 2D cubic CGL equation~\eqref{eq:ccgl} with periodic
      boundary conditions and different integrators. The total number of 
      degrees of freedom is $N=256^2$. In the table, the cross symbol $\times$ 
      means NaN/Inf value returned.
      See also Figure~\ref{fig:cgl_cubic_2D_fourier_cpudiag} for a graphical representation.}
  \label{tab:cgl_cubic_2D_fourier}
\end{table}
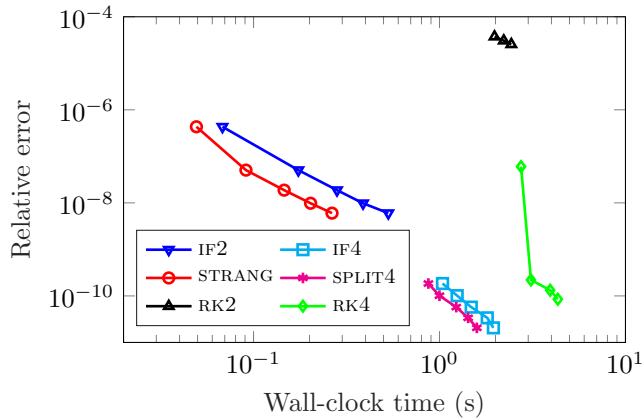
\begin{figure}
  \centering
%
%
%
\begin{tikzpicture}

\begin{axis}[%
width=2.6in,
height=1.7in,
at={(0.758in,0.481in)},
scale only axis,
xmode=log,
xmin=0.02,
xmax=10,
xminorticks=true,
xlabel style={font=\color{white!15!black}},
ymode=log,
ymin=1e-11,
ymax=1e-04,
yminorticks=true,
xlabel = {Wall-clock time (s)},
ylabel = {Relative error},
ylabel style={font=\color{white!15!black}},
axis background/.style={fill=white},
legend style={at={(0.025,0.05)}, anchor=south west, legend cell align=left, align=left, draw=white!15!black, font=\footnotesize},
legend columns = 2
]
\addplot [color=blue, line width=1pt, mark size = 2pt,mark=triangle, mark options={solid, rotate=180, blue}]
  table[row sep=crcr]{%
0.068061	4.32241290873405e-07\\
0.174179	5.05525224828908e-08\\
0.281071000000001	1.87657359967463e-08\\
0.388111999999999	9.80884751523611e-09\\
0.530536000000001	6.0141842877127e-09\\
};
\addlegendentry{\textsc{if2}}

\addplot [color=cyan, line width=1pt, mark size = 2pt,mark=square, mark options={solid, cyan}]
  table[row sep=crcr]{%
1.037871	1.84637781150687e-10\\
1.244201	1.01216330063991e-10\\
1.481196	5.74169934032493e-11\\
1.806753	3.37807352318231e-11\\
1.949446	2.05949798131755e-11\\
};
\addlegendentry{\textsc{if4}}

\addplot [color=red, line width=1pt, mark size = 2pt,mark=o, mark options={solid, red}]
  table[row sep=crcr]{%
0.049299	4.32231907325279e-07\\
0.0910429999999998	5.05521863431542e-08\\
0.146123	1.87656677924465e-08\\
0.202713	9.80882470590973e-09\\
0.264288	6.01417509186981e-09\\
};
\addlegendentry{\textsc{strang}}

\addplot [color=magenta,line width=1pt, mark size = 2pt, mark=asterisk, mark options={solid, magenta}]
  table[row sep=crcr]{%
0.871719000000002	1.84625290961171e-10\\
0.999693	1.01193213340016e-10\\
1.230777	5.74150604213386e-11\\
1.425227	3.37843184331536e-11\\
1.588141	2.06001141498515e-11\\
};
\addlegendentry{\textsc{split4}}

\addplot [color=black, line width=1pt, mark size = 2pt,mark=triangle, mark options={solid, black}]
  table[row sep=crcr]{%
1.9763e+00	3.7324e-05\\
2.2167e+00	3.0596e-05\\
2.4338e+00	2.5536e-05\\
};
\addlegendentry{\textsc{rk2}}

\addplot [color=green, line width=1pt, mark size = 2pt,mark=diamond, mark options={solid, green}]
  table[row sep=crcr]{%
2.74508	6.01494354373691e-08\\
3.10362400000001	2.20504016242741e-10\\
3.937045	1.32648891947163e-10\\
4.328766	8.49063200548815e-11\\
};
\addlegendentry{\textsc{rk4}}

\end{axis}

\end{tikzpicture}%
  \caption{Results for the simulation of the 2D cubic CGL
    equation~\eqref{eq:ccgl} with periodic
    boundary conditions and
  $N=256^2$ spatial discretization points. The number
  of time steps for each integrator is reported in
  Table~\ref{tab:cgl_cubic_2D_fourier}. The final simulation time is $T=6$.}
  \label{fig:cgl_cubic_2D_fourier_cpudiag}
\end{figure}

Again, for the first experiment we set the final simulation time to $T=6$. The
outcome is reported in Table~\ref{tab:cgl_cubic_2D_fourier} and graphically 
presented in Figure~\ref{fig:cgl_cubic_2D_fourier_cpudiag}.
Similar considerations as for the previous example can be drawn.
{\color{black}The explicit Runge--Kutta integrators perform poorly.}
The second-order exponential-type methods reach comparable accuracies using
the same number of time steps, and similarly for the fourth-order ones.
{\color{black}
Splitting methods are slightly more efficient in terms of computational time.
Moreover, the advantage of using fourth-order exponential-type
methods for stringent accuracies is even more evident and,
for a wide range of tolerances, the best integrator
is again \textsc{split4}.}

\begin{figure}[!htb]
  \centering
  \input{cgl_cubic_2D_fourier_dyn.tex}
  \caption{Dynamics of $\lvert u \rvert$ for the 2D cubic CGL
    equation~\eqref{eq:ccgl} with periodic
    boundary conditions and
  $N=256^2$ spatial discretization points. The time integrator is \textsc{split4}.
  The final times are $T=30$ (top left), $T=60$ (top right),
  $T=120$ (bottom left), and $T=240$ (bottom right) with a number of time steps
  equal to
  $m=50$, $m=100$, $m=200$, and $m=400$, respectively.
  The wall-clock times
  of the simulations are $0.33$, $0.59$, $1.15$, and $2.20$ seconds,
  respectively.}
  \label{fig:cgl_cubic_2D_fourier_dyn}
\end{figure}

For the second experiment, we show in Figure~\ref{fig:cgl_cubic_2D_fourier_dyn}
the dynamics of $\lvert u \rvert$ at final times $T=30$, $T=60$, $T=120$, and $T=240$.
The time integrator is \textsc{split4} with a number of time steps equal to
$m=50$, $m=100$, $m=200$, and $m=400$, respectively.
Also in this case, the behavior matches with the theoretical expectations. In
particular, the dynamics is deeply different from the case of homogeneous Dirichlet
boundary conditions. Indeed, the initially formed defects move and merge but
now cross the boundaries in a periodic fashion, and no frozen state is observed.
In terms of wall-clock time, the simulation with the
largest number of times steps
executed in about two seconds.
\subsection{3D cubic CGL equation with homogeneous Dirichlet--Neumann boundary conditions}
\label{sec:CCGL3D}
We now consider the cubic CGL equation
\begin{equation}\label{eq:ccgl3D}
  \partial_tu=(\alpha_1+\rmi\beta_1)\Delta u+
  \alpha_2 u+
  (\alpha_3+\rmi\beta_3)\lvert u\rvert^2u,
\end{equation}
in the three-dimensional domain {\color{black}$\Omega=(0,100)^3$} with parameters
$\alpha_1=1$, $\beta_1=2$, $\alpha_2=1$, $\alpha_3=-1$,
and $\beta_3=0.2$. The number of discretization points is equal to 
$n=128$ in all the three directions (total number of degrees of freedom
$N=128^3$).
As for the two-dimensional examples of Section~\ref{sec:CCGL2D}, the initial
datum is a standard normally distributed random obtained as 
\verb+randn(n,n,n)/5000+ with \verb+rng('default')+. Also, for the splitting
methods the nonlinear flow can be again analytically computed as in
formula~\eqref{eq:exsolcub}.
The equation is completed with homogeneous Dirichlet boundary conditions in the
portion of boundary given by $\{(x_1,x_2,x_3)\colon x_1x_2x_3=0\}$, 
and homogeneous Neumann elsewhere (see Reference~\cite{EHGP01} for a similar
example in two spatial dimensions).
The discretization in space is then performed
with fourth-order finite differences using
semi-internal nodes in each direction.
The unidirectional discretization matrix is in this case
\begin{equation*}
  \frac{1}{12h^2}\begin{bmatrix}
    -15 & -4  & 14 & -6 & 1 & 0 & \ldots & \ldots & 0\\
    16 & -30 & 16 & -1 & 0 & \ldots & \ldots & \ldots & 0\\
    -1 & 16 & -30 & 16 & -1 & 0 & \ldots & \ldots & 0\\
    0 & -1 & 16 & -30 & 16 & -1 & 0 & \ldots & 0\\
    \vdots & \ddots & \ddots & \ddots & \ddots & \ddots &\ddots & \ddots & \vdots\\
    0 & \ldots & 0 & -1 & 16 & -30 & 16 & -1 & 0\\
    0 & \ldots & \ldots & 0 & -1 & 16 & -30 & 16 & -1\\
    0 & \ldots & 0 & 1 & -6 & 14 & -4 & -15 & 10\\
    0 & \ldots & \ldots & 0 & 1 & -\frac{8}{3} & -6 & 56 & -\frac{145}{3}
  \end{bmatrix},
\end{equation*}
in all directions, where $h=100/n$ is the spatial step size. 
The first row (which corresponds to the first internal point of the discretized 
domain) is the same as for the one in Section~\ref{sec:CCGL2Dhd}, while the last
two (which corresponds to the last two points of the discretized domain) 
have been suitably modified to take into account a fourth-order discretization
of the homogeneous Neumann boundary conditions and the lack of enough 
surrounding points for the usual centered stencil.


\begin{table}[!htb]
  \centering
  {\small
  \begin{tabular}{llll|llll}
    \hline
    \multicolumn{4}{l|}{\textsc{if2}} &                     \multicolumn{4}{l}{\textsc{if4}} \\                     
        steps & time (s) & error & order &             steps & time (s) & error & order \\           
    \hline                                             
    25 & $3.964\rme{+00}$ & $1.235\rme{-03}$ & ---  &    50 & $2.247\rme{+01}$ & $1.889\rme{-07}$ & --- \\
    75 & $1.175\rme{+01}$ & $1.349\rme{-04}$ & 2.015 &   100 & $4.506\rme{+01}$ & $1.143\rme{-08}$ & 4.047 \\ 
    125 & $1.956\rme{+01}$ & $4.834\rme{-05}$ & 2.009 &  150 & $6.879\rme{+01}$ & $2.194\rme{-09}$ & 4.071 \\
    175 & $2.812\rme{+01}$ & $2.461\rme{-05}$ & 2.007 &  200 & $9.175\rme{+01}$ & $6.868\rme{-10}$ & 4.036 \\
    225 & $3.683\rme{+01}$ & $1.4869\rme{-05}$ & 2.005 & 250 & $1.206\rme{+02}$ & $2.803\rme{-10}$ & 4.017 \\
    \hline
    \hline
 \multicolumn{4}{l|}{\textsc{strang}}                      & \multicolumn{4}{l}{\textsc{split4}}\\
\hline
steps & time (s) & error & order                      & steps & time (s) & error & order\\
\hline                                                       
 25 & $2.342\rme{+00}$ & $1.190\rme{-03}$ & ---          & 50 & $1.852\rme{+01}$& $1.848\rme{-07}$ & ---     \\
 75 & $6.855\rme{+00}$ & $1.335\rme{-04}$ & 1.991        & 100 & $3.237\rme{+01}$ & $1.114\rme{-08}$ & 4.052 \\
 125 & $1.110\rme{+01}$ & $4.809\rme{-05}$ & 1.999      & 150 & $4.714\rme{+01}$ & $2.135\rme{-09}$ & 4.075 \\
 175 & $1.567\rme{+01}$ & $2.454\rme{-05}$ & 1.999       & 200 & $6.261\rme{+01}$ & $6.679\rme{-10}$ & 4.039 \\
 225 & $1.993\rme{+01}$ & $1.485\rme{-05}$ & 2.000       & 250 & $7.815\rme{+01}$ & $2.724\rme{-10}$ & 4.019 \\
    \hline
    \hline
    \multicolumn{4}{l|}{\textsc{rk2}}&
    \multicolumn{4}{l}{\textsc{rk4}}\\
    steps & time (s) & error & order &
    steps & time (s) & error & order \\
    \hline
    115 & $\times$ & $\times$  & --- &
    115 & $\times$ & $\times$ & --- \\
    215 & $\times$ & $\times$ & $\times$ &
    215 & $8.885\rme{+01}$ & $8.108\rme{-01}$ & $\times$\\
    315 & $6.137\rme{+01}$ & $1.292\rme{-03}$  & $\times$ &
    315 & $1.311\rme{+02}$ & $5.817\rme{-08}$ & 43.07\\
    415 & $8.163\rme{+01}$&  $7.484\rme{-04}$ & 1.979 &
    415 & $1.700\rme{+02}$ & $1.942\rme{-08}$ & 3.980\\
    515 & $1.001\rme{+02}$ & $4.876\rme{-04}$ & 1.984 &
    515 & $2.108\rme{+02}$ & $8.214\rme{-09}$ & 3.984\\
    \hline
  \end{tabular}%
  }
  \caption{Number of time steps, wall-clock time (in seconds),
  relative error at final time $T=10$, and observed numerical order of convergence
      for the solution of the 3D cubic CGL equation~\eqref{eq:ccgl3D} with homogeneous
      Dirichlet--Neumann boundary conditions and different integrators. 
      The total number of degrees of freedom is $N=128^3$. In the table, the cross symbol $\times$ 
      means NaN/Inf value returned.
      See also Figure~\ref{fig:cgl_cubic_3D_FD_cpudiag} for a graphical representation.}
  \label{tab:cgl_cubic_3D_FD}
\end{table}
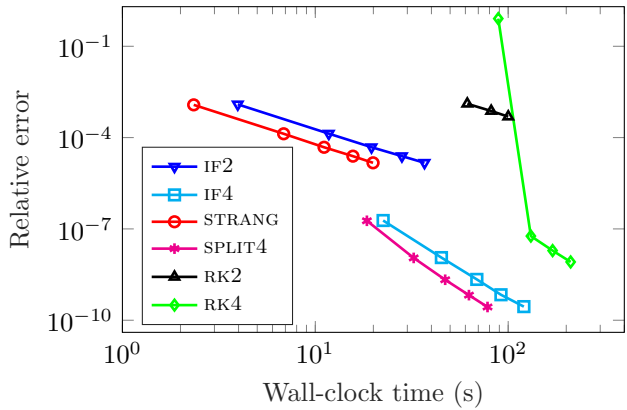
\begin{figure}[!htb]
  \centering
%
%
%
\begin{tikzpicture}

\begin{axis}[%
width=2.6in,
height=1.7in,
at={(0.769in,0.477in)},
scale only axis,
xmode=log,
xmin=1,
xmax=400,
xminorticks=true,
xlabel style={font=\color{white!15!black}},
ymode=log,
ymin=4e-11,
ymax=2,
yminorticks=true,
ylabel style={font=\color{white!15!black}},
xlabel = {Wall-clock time (s)},
ylabel = {Relative error},
axis background/.style={fill=white},
legend style={at={(0.34,0.57)}, legend cell align=left, align=left, draw=white!15!black, font=\footnotesize}
]
\addplot [color=blue, line width=1pt, mark size = 2pt, mark=triangle, mark options={solid, rotate=180, blue}]
  table[row sep=crcr]{%
3.963898	0.00123469468823523\\
11.746858	0.000134923939682688\\
19.563847	4.83411928529825e-05\\
28.122871	2.46103727523173e-05\\
36.828305	1.48693453718449e-05\\
};
\addlegendentry{\textsc{if2}}

\addplot [color=cyan, line width=1pt, mark size = 2pt,mark=square, mark options={solid, cyan}]
  table[row sep=crcr]{%
22.472502	1.88898904103839e-07\\
45.058831	1.14289389483707e-08\\
68.7861490000001	2.19348573110808e-09\\
91.7525540000001	6.86833894297039e-10\\
120.578171	2.80275078387375e-10\\
};
\addlegendentry{\textsc{if4}}

\addplot [color=red, line width=1pt, mark size = 2pt,mark=o, mark options={solid, red}]
  table[row sep=crcr]{%
2.342153	0.00118974338268952\\
6.854754	0.000133479074748457\\
11.104347	4.80874790434939e-05\\
15.673423	2.45392666520792e-05\\
19.92648	1.4845928628837e-05\\
};
\addlegendentry{\textsc{strang}}

\addplot [color=magenta, line width=1pt, mark size = 2pt,mark=asterisk, mark options={solid, magenta}]
  table[row sep=crcr]{%
18.515622	1.84826668013761e-07\\
32.3684900000001	1.11424468956564e-08\\
47.139124	2.13476816418379e-09\\
62.607011	6.67866881121883e-10\\
78.1542080000001	2.7237597436364e-10\\
};
\addlegendentry{\textsc{split4}}

\addplot [color=black, line width=1pt, mark size = 2pt,mark=triangle, mark options={solid, black}]
  table[row sep=crcr]{%
6.1373e+01	1.2916e-03\\
8.1630e+01	7.4837e-04\\
1.0008e+02	4.8764e-04\\
};
\addlegendentry{\textsc{rk2}}

\addplot [color=green, line width=1pt, mark size = 2pt,mark=diamond, mark options={solid, green}]
  table[row sep=crcr]{%
88.849373	8.107987277455518e-01\\
131.094346	5.817404981526199e-08\\
169.950732	1.941601312656044e-08\\
210.790786	8.214490377625786e-09\\
};
\addlegendentry{\textsc{rk4}}

\end{axis}

\end{tikzpicture}%
  \caption{Results for the simulation of the 3D cubic CGL
    equation~\eqref{eq:ccgl3D} with homogeneous Dirichlet--Neumann
    boundary conditions and
    $N=128^3$ spatial discretization points. The number
    of time steps for each integrator is reported in
    Table~\ref{tab:cgl_cubic_3D_FD}. The final simulation time is $T=10$.}
  \label{fig:cgl_cubic_3D_FD_cpudiag}
\end{figure}

For the first experiment we set the final time to $T=10$ and test the 
performances of the different integrators. The results are summarized 
in Table~\ref{tab:cgl_cubic_3D_FD} and depicted in
Figure~\ref{fig:cgl_cubic_3D_FD_cpudiag}.
The conclusions are equivalent to the ones drawn for the numerical examples
in Section~\ref{sec:CCGL2D}. {\color{black}In particular,
  the explicit Runge--Kutta methods
suffer from a time step size restriction due to the stiffness,
and
a slight advantage for the splitting methods compared
to the Lawson schemes is observed.}

\begin{figure}[!htb]
  \centering
  \input{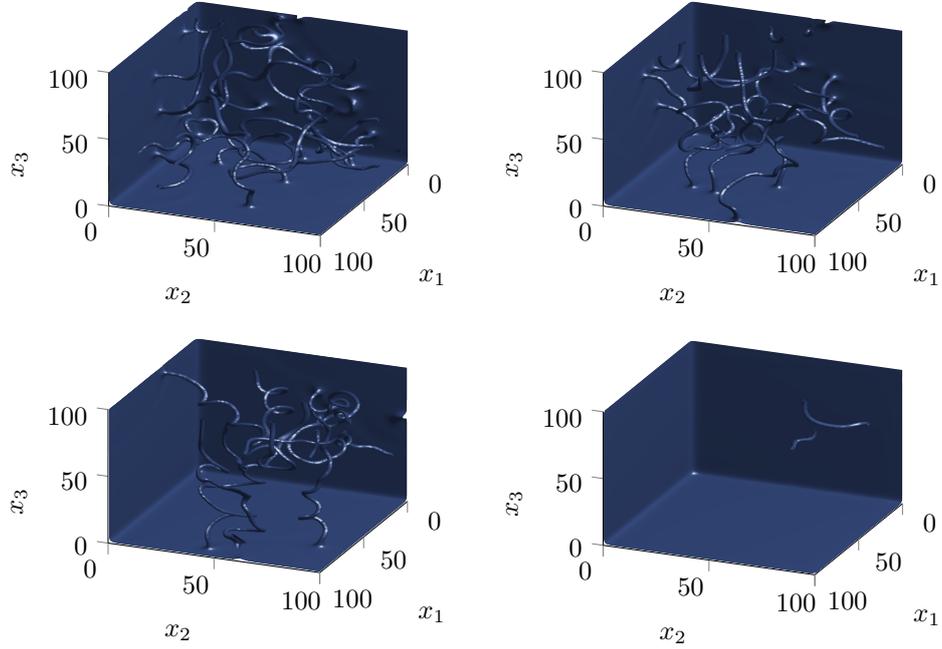}
  \caption{Dynamics of $\lvert u \rvert$ for the 3D cubic CGL
    equation~\eqref{eq:ccgl3D} with homogeneous Dirichlet--Neumann
    boundary conditions and
  $N=128^3$ spatial discretization points (isosurface level $0.3$).
  The time integrator is \textsc{split4}. 
  The final times are $T=100$ (top left), $T=200$ (top right), 
  $T=300$ (bottom left), and $T=400$ (bottom right) with a number of time steps 
  equal to
  $m=1000$, $m=2000$, $m=3000$, and $m=4000$, respectively. 
  The wall-clock times
  of the simulations are $306.82$, $630.23$, $917.55$, and $1272.67$ seconds,
  respectively.}
  \label{fig:cgl_cubic_3D_FD_dyn}
\end{figure}

Concerning the second experiment, we investigate the dynamics of 
$\lvert u \rvert$ for the three-dimensional equation~\eqref{eq:ccgl3D} using
\textsc{split4} as
time integrator. To this aim we perform simulations up to
the final times $T=100$, $T=200$, $T=300$, and $T=400$ with a number of time
steps equal to $m=1000$, $m=2000$, $m=3000$, and $m=4000$, respectively.
We show in Figure~\ref{fig:cgl_cubic_3D_FD_dyn} the corresponding 
isosurfaces at level $0.3$. As we expect, in the three-dimensional scenario the
defects have the shape of vortex filaments (see, for instance,
References~\cite{ABK98,RCK08}),
which move and evolve as the time
proceeds. More specifically, they are pushed away from the faces of the cube 
in which homogeneous Dirichlet boundary conditions are prescribed, and they are 
ejected from the opposite faces (in which homogeneous Neumann boundary 
conditions are set). In fact, if we let evolve the system after time $T=400$,
no more defects are present in the considered domain. This is in line with what
reported for the two-dimensional case in the literature (see, for instance,
Reference~\cite{EHGP01}).

\subsection{3D cubic-quintic CGL equation with periodic boundary conditions}
\begin{figure}[!htb]
  \centering
  \input{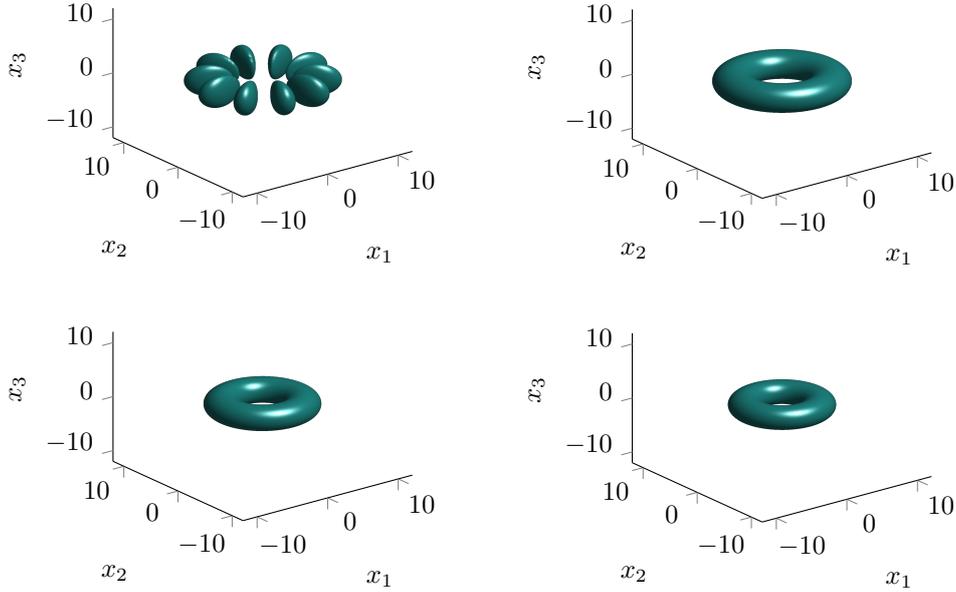}
  \caption{Dynamics of $\lvert u \rvert^2$ for the 3D
  cubic-quintic CGL
    equation~\eqref{eq:cqgl} with periodic
    boundary conditions and
    $N=128^3$ spatial discretization points
    (isosurface level 0.5). The time integrator is \textsc{if4}. 
    We depict the initial condition (top left plot), and the isosurfaces at 
    final times $T=20$ (top right plot), $T=40$ (bottom left plot), and $T=80$
    (bottom right plot). The number of time steps is equal to 
  $m=80$, $m=160$, $m=320$, respectively. The wall-clock times
  of the simulations are $38.93$, $81.38$, and $159.65$ seconds,
  respectively.}
  \label{fig:cgl_cubicquintic_3D_fourier_dyn}
\end{figure}

We now turn our attention to the following cubic-quintic CGL equation
\begin{equation}\label{eq:cqgl}
  \partial_tu=(\alpha_1+\rmi\beta_1)\Delta u+
  \alpha_2 u+
  (\alpha_3+\rmi\beta_3)\lvert u\rvert^2u+
    (\alpha_4+\rmi\beta_4)\lvert u\rvert^4u,
\end{equation}
in the cube {\color{black}$\Omega=(-12,12)^3$} completed with periodic boundary
conditions (see References~\cite{LHQW09,ZYW16}). The parameters are set to
$\alpha_1=1/2$, $\beta_1=1/2$, $\alpha_2=-1/2$, $\alpha_3=2.52$,
$\beta_3=1$, $\alpha_4=-1$, and $\beta_4=-0.11$.
Let us define $\rho=\sqrt{x_1^2+x_2^2}$ and
$\theta=\at2(x_2,x_1)$. Then, we take as initial datum
\begin{equation*}
u_0=  \delta\sech\left(\frac{\sqrt{(\rho-\rho_0)^2+x_3^2}}{\omega}\right)
  \cos(\eta\theta)\exp(\rmi \kappa\theta),
\end{equation*}
with $\delta=1.2$, $\rho_0=6$, $\omega=2.5$, $\eta=5$, and $\kappa=3$. This
initial condition corresponds to a \textit{necklace-ring} with 10 \textit{beads},
see the top left plot in Figure~\ref{fig:cgl_cubicquintic_3D_fourier_dyn} 
for a graphical representation. As space discretization we employ a Fourier
pseudospectral approach with $n=128$ modes along each direction, obtaining
a total number of degrees of freedom equal to $N=128^3$.

\begin{table}[!htb]
  \centering
  {\small
  \begin{tabular}{llll|llll}
    \hline
    \multicolumn{4}{l|}{\textsc{if2}} &                        \multicolumn{4}{l}{\textsc{if4}} \\                      
        steps & time (s) & error & order &                    steps & time (s) & error & order \\             
    \hline                                                
    1000 & $1.976\rme{+02}$ & $2.197\rme{-04}$ & ---  &   300 & $1.417\rme{+02}$ & $8.617\rme{-07}$ & ---  \\ 
    1500 & $2.952\rme{+02}$ & $9.859\rme{-05}$ & 1.977 &  450 & $2.130\rme{+02}$ & $1.659\rme{-07}$ & 4.064 \\
    2000 & $3.940\rme{+02}$ & $5.572\rme{-05}$ & 1.984 &  600 & $2.842\rme{+02}$ & $5.186\rme{-08}$ & 4.042 \\
    2500 & $4.956\rme{+02}$ & $3.576\rme{-05}$ & 1.987 &  750 & $3.554\rme{+02}$ & $2.109\rme{-08}$ & 4.031 \\
    3000 & $5.935\rme{+02}$ & $2.488\rme{-05}$ & 1.990 &  900 & $4.272\rme{+02}$ & $1.013\rme{-08}$ & 4.025 \\
    \hline
    \hline
    \multicolumn{4}{l|}{\textsc{strang}} &                            \multicolumn{4}{l}{\textsc{split4}} \\                    
    steps & time (s) & error & order &                      steps & time (s) & error & order \\                 
    \hline                                                                                                     
    1000 & $2.080\rme{+02}$ & $9.212\rme{-05}$ & ---   &    300 & $2.804\rme{+02}$ & $3.005\rme{-07}$ & ---   \\
    1500 & $3.126\rme{+02}$ & $4.094\rme{-05}$ & 2.000 &    450 & $4.213\rme{+02}$ & $5.942\rme{-08}$ & 3.998 \\
    2000 & $4.173\rme{+02}$ & $2.303\rme{-05}$ & 2.000 &    600 & $5.623\rme{+02}$ & $1.878\rme{-08}$ & 4.004 \\
    2500 & $5.226\rme{+02}$ & $1.474\rme{-05}$ & 2.000 &    750 & $7.020\rme{+02}$ & $7.672\rme{-09}$ & 4.011 \\
    3000 & $6.284\rme{+02}$ & $1.023\rme{-05}$ & 2.000 &    900 & $8.426\rme{+02}$ & $3.685\rme{-09}$ & 4.023 \\
    \hline
    \hline
    \multicolumn{4}{l|}{\textsc{strang\_3t}} &                        \multicolumn{4}{l}{\textsc{split4\_3t}}\\                  
    steps & time(s) & error & order &                     steps & time(s) & error & order\\                  
    \hline                                                                                                   
    1000 & $1.237\rme{+02}$ & $6.161\rme{-04}$ & --- &    300 & $1.130\rme{+02}$ & $3.189\rme{-05}$ & ---\\  
    1500 & $1.706\rme{+02}$ & $2.738\rme{-04}$ & 2.000 &  450 & $1.691\rme{+02}$ & $6.378\rme{-06}$ & 3.969\\
    2000 & $2.268\rme{+02}$ & $1.540\rme{-04}$ & 2.000 &  600 & $2.255\rme{+02}$ & $2.034\rme{-06}$ & 3.973\\
    2500 & $2.830\rme{+02}$ & $9.858\rme{-05}$ & 2.000 &  750 & $2.830\rme{+02}$ & $8.375\rme{-07}$ & 3.977\\
    3000 & $3.396\rme{+02}$ & $6.846\rme{-05}$ & 2.000 &  900 & $3.379\rme{+02}$ & $4.054\rme{-07}$ & 3.980\\
    \hline
    \hline
    \multicolumn{4}{l|}{\textsc{rk2}} &
    \multicolumn{4}{l}{\textsc{rk4}}\\
    steps & time (s) & error & order &
    steps & time (s) & error & order \\
    \hline
    1200 & $\times$ & $\times$ & --- &
    1000 & $\times$ & $\times$ & --- \\
    1300 & $\times$  & $\times$  & $\times$  &
    1100 & $4.492\rme{+02}$ & $1.273\rme{-07}$ & $\times$\\
    1400 & $2.918\rme{+02}$  & $1.386\rme{-05}$ & $\times$ &
    1200 & $4.901\rme{+02}$ & $5.698\rme{-10}$ & 62.17\\
    1500 & $3.127\rme{+02}$ & $1.208\rme{-05}$ & 2.002 &
    1300 & $5.320\rme{+02}$ & $4.133\rme{-10}$ & 4.012\\
    1600 & $3.347\rme{+02}$ & $1.061\rme{-05}$ & 2.002 &
    1400 & $5.753\rme{+02}$ & $3.077\rme{-10}$ & 3.982\\
    \hline
  \end{tabular}%
  }
  \caption{Number of time steps, wall-clock time (in seconds),
  relative error at final time $T=5$, and observed numerical order of convergence
      for the solution of the 3D cubic-quintic CGL equation~\eqref{eq:cqgl} with periodic
      boundary conditions and different integrators. The total number of 
      degrees of freedom is $N=128^3$. In the table, the cross symbol $\times$ 
      means NaN/Inf value returned.
      See also Figure~\ref{fig:cgl_cubicquintic_3D_fourier_cpudiag} for a graphical representation.}
  \label{tab:cgl_cubicquintic_3D_fourier}
\end{table}
\begin{figure}[!htb]
  \centering
%
%
%
\begin{tikzpicture}

\begin{axis}[%
width=2.8in,
height=1.7in,
at={(0.769in,0.477in)},
scale only axis,
xmode=log,
xmin=99,
xmax=1001,
xminorticks=true,
xlabel style={font=\color{white!15!black}},
ymode=log,
ymin=2e-12,
ymax=1e-3,
yminorticks=true,
ylabel style={font=\color{white!15!black}},
xlabel = {Wall-clock time (s)},
ylabel = {Relative error},
axis background/.style={fill=white},
legend style={at={(0.65,0.38)}, legend cell align=left, align=left, draw=white!15!black, font=\scriptsize},
legend columns=2
]
\addplot [color=blue, line width=1pt, mark size = 2pt, mark=triangle, mark options={solid, rotate=180, blue}]
  table[row sep=crcr]{%
197.607392	0.000219740408156955\\
295.185331	9.85891281640069e-05\\
393.966763	5.57184399408641e-05\\
495.556978	3.57607094536931e-05\\
593.505501000001	2.48806199039194e-05\\
};
\addlegendentry{\textsc{if2}}

\addplot [color=cyan, line width=1pt, mark size = 2pt,mark=square, mark options={solid, cyan}]
  table[row sep=crcr]{%
141.688032	8.61721229897972e-07\\
212.960896	1.65888683922702e-07\\
284.208978	5.18557001728487e-08\\
355.408428	2.10916902728832e-08\\
427.170238	1.01262090939577e-08\\
};
\addlegendentry{\textsc{if4}}

\addplot [color=red, line width=1pt, mark size = 2pt,mark=o, mark options={solid, red}]
  table[row sep=crcr]{%
208.025433	9.21216568080117e-05\\
312.591374	4.09393365309816e-05\\
417.30808	2.30276602456665e-05\\
522.586266	1.47374847602681e-05\\
628.397736	1.02342778333317e-05\\
};
\addlegendentry{\textsc{strang}}

\addplot [color=magenta, line width=1pt, mark size = 2pt,mark=asterisk, mark options={solid, magenta}]
  table[row sep=crcr]{%
280.352352	3.00505104935729e-07\\
421.267841	5.94187494152577e-08\\
562.306558999999	1.87767241268078e-08\\
702.006030000001	7.67222280888924e-09\\
842.639863	3.68451549596402e-09\\
};
\addlegendentry{\textsc{split4}}

\addplot [color=teal, line width=1pt, mark size = 2pt,mark=+, mark options={solid, teal}]
  table[row sep=crcr]{%
123.717538	0.000616128568698997\\
170.577269	0.000273827218917739\\
226.768314	0.000154026307029259\\
283.001187	9.8576396925393e-05\\
339.611885	6.84556692758333e-05\\
};
\addlegendentry{\textsc{strang\_3t}}
\addplot [color=violet, line width=1pt, mark size = 2pt,mark=star, mark options={solid, violet}]
  table[row sep=crcr]{%
113.000543	3.18898486805101e-05\\
169.062043	6.37849245070653e-06\\
225.450657	2.03402591042135e-06\\
283.034949	8.37486654793441e-07\\
337.888056	4.05369229338754e-07\\
};
\addlegendentry{\textsc{split4\_3t}}

\addplot [color=black, line width=1pt, mark size = 2pt,mark=triangle, mark options={solid, black}]
  table[row sep=crcr]{%
2.9180e+02	1.3864e-05\\
3.1266e+02	1.2076e-05\\
3.3472e+02	1.0612e-05\\
};
\addlegendentry{\textsc{rk2}}

\addplot [color=green, line width=1pt, mark size = 2pt,mark=diamond, mark options={solid, green}]
  table[row sep=crcr]{%
449.153162	1.27335643976485e-07\\
490.115522000001	5.69751780212045e-10\\
531.98705	4.13258706799146e-10\\
575.258136000001	3.07659706185418e-10\\
};
\addlegendentry{\textsc{rk4}}

\end{axis}

\end{tikzpicture}%
  \caption{Results for the simulation of the 3D cubic-quintic CGL
    equation~\eqref{eq:cqgl} with periodic
    boundary conditions and
    $N=128^3$ spatial discretization points. The number
    of time steps for each integrator is reported in
    Table~\ref{tab:cgl_cubicquintic_3D_fourier}. The final simulation time is $T=5$.}
  \label{fig:cgl_cubicquintic_3D_fourier_cpudiag}
\end{figure}
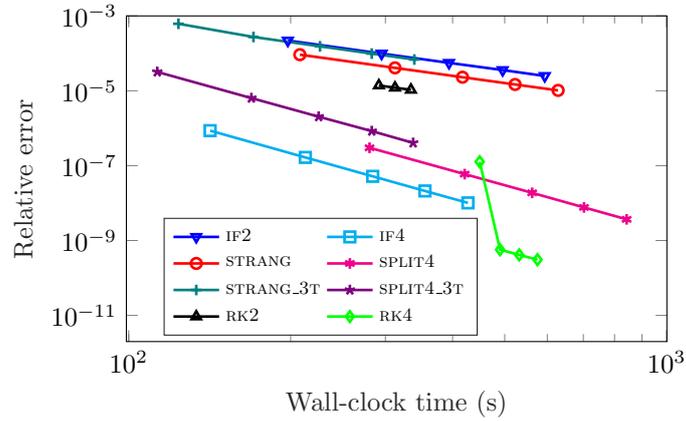

For the time integration in the first experiment, we consider all the 
marching schemes with final time set to $T=5$.
Concerning specifically the splitting methods, notice that contrary to 
the cubic nonlinearity case considered in
Sections~\ref{sec:CCGL2D} and~\ref{sec:CCGL3D}, here there is no exact solution
for the nonlinear flow (see also the discussion in Reference~\cite{WZ13}).
In fact, in the time marching we numerically compute an approximation 
of it with a single step of \textsc{rk4}, since there is no source of stiffness 
and it overall gave the best results (the linear part is obviously integrated
exactly). An alternative approach is to separate the nonlinearity into a cubic
part and a quintic one, of which analytical solutions are known. Indeed, for the
former we have the exact solution given by formula~\eqref{eq:exsolcub}, while
for the latter if we consider
\begin{equation*}
  \check{\bu}'(t) = (\alpha_4+\rmi\beta_4)\lvert \check{\bu}(t)\rvert^4\check{\bu}(t)
\end{equation*}
with initial datum $\check{\bu}_0$ we have the exact solution
\begin{equation*}
  \check{\bu}(t) = \phi_{t}^{\mathrm q}(\check{\bu}_0)=\rme^{-\frac{\alpha_4+\rmi\beta_4}{4\alpha_4}
                    \log\lvert1-4\alpha_4\lvert \check{\bu}_0 \rvert^4t\rvert}
                    \check{\bu}_0.
\end{equation*}
Embedding this procedure in a three-term splitting approach gives rise in principle
to valuable alternatives for the time marching. 
In particular, for the second-order Strang splitting we consider
\begin{equation*}
  \bu^{n+1} = 
  \mathcal{S}^3_\tau(\bu_n) =
  \phi^{\mathrm{q}}_{\tau/2}(\phi^{\mathrm{c}}_{\tau/2}(\rme^{\tau K}\phi^{\mathrm{c}}_{\tau/2}(\phi^{\mathrm{q}}_{\tau/2}(\bu_n)))),
\end{equation*}
where $\phi^{\mathrm{c}}$ and $\phi^{\mathrm{q}}$ denote the exact flows for
the cubic and the quintic parts, respectively. The scheme is labeled as
\textsc{strang\_3t}. Then, similarly to formula~\eqref{eq:richextrap}, 
if we perform the Richardson extrapolation we obtain the fourth-order method
\begin{equation*}
\bu_{n+1} = \frac{4}{3}\mathcal{S}^3_{\tau/2}(\mathcal{S}^3_{\tau/2}(\bu_n))-
  \frac{1}{3}\mathcal{S}^3_{\tau}(\bu_n),
\end{equation*}
that we label \textsc{split4\_3t}.
The outcome of the experiments is summarized in
Table~\ref{tab:cgl_cubicquintic_3D_fourier} and plotted in 
Figure~\ref{fig:cgl_cubicquintic_3D_fourier_cpudiag}.
As for all the previous examples, we notice that all the time integrators show
the expected order of convergence. Also, the {\color{black} explicit Runge--Kutta
methods again show} a severe time step size restriction, which 
{\color{black}makes them} not of practical usage
for fast simulations with a low number of time
steps and large number of degrees of freedom.
{\color{black}Concerning the second-order exponential-type methods,
  we observe from the table that the three-term 
splitting \textsc{strang\_3t} is cheaper than the \textsc{strang} method
with nonlinear flux computed numerically, as expected. On the other hand, the 
relative error obtained is almost seven times larger. This makes
the \textsc{strang} scheme the preferred one in the precision-work diagram,
since the \textsc{if2} method performs similarly to \textsc{strang\_3t}.}
For the exponential-type methods of fourth order, we observe that
the three-term splitting approach \textsc{split4\_3t} does not pay off compared
to \textsc{split4}, since the lower computational cost is surpassed by a 
considerable higher error obtained. As a matter of fact, the \textsc{if4}
is more efficient than {\color{black}all the other exponential-type methods.}
We then conclude that for the example under consideration the integrator which
performs best is \textsc{if4}. 

In the second experiment we present the dynamics of the model for different final
times. To this aim, we select as time integrator \textsc{if4}, and we consider
the final simulation times $T=20$, $T=40$, and $T=80$ with a number of time steps
equal to $m=80$, $m=160$, and $m=320$, respectively.
The results, depicted in Figure~\ref{fig:cgl_cubicquintic_3D_fourier_dyn},
  are coherent with what already presented in the literature, see 
References~\cite{LHQW09,ZYW16}. In particular, the magnitude of the coefficients
in the equation makes evolve the initial necklace-ring to a torus, which shrinks as
time evolves. After time $T=120$ the isosurface does not show qualitatively
modifications, suggesting that the stationary state has been reached.
As already observed previously in the numerical
experiments, we notice that also here the
simulation wall-clock time of a single step is the same for each simulation,
since the underlying marching method is a direct algorithm. The simulation
with the largest number of time steps takes about 160 seconds.

\subsection{Coupled 2D cubic-quintic CGL equation with periodic boundary conditions}
We consider as last example the interaction of two quasi-1D dissipative
solitons modeled by a system of two coupled cubic-quintic CGL equations
(see Reference~\cite{DB14})
\begin{equation}\label{eq:coupled}
  \left\{
\begin{aligned}
  \partial_tu=&\,\alpha_0\partial_{x_1} u+(\alpha_1+\rmi\beta_1)\Delta u+
  \alpha_2u\\
  &+(\alpha_3+\rmi\beta_3)\lvert u\rvert^2u+
  (\alpha_4+\rmi\beta_4)\lvert u\rvert^4u+
  \alpha_5\lvert v\rvert^2u,\\
  \partial_tv=&-\alpha_0\partial_{x_1} v+(\alpha_1+\rmi\beta_1)\Delta v+
  \alpha_2 v\\
  &+(\alpha_3+\rmi\beta_3)\lvert v\rvert^2v+
  (\alpha_4+\rmi\beta_4)\lvert v\rvert^4v+
  \alpha_5\lvert u\rvert^2v.
\end{aligned}\right.
\end{equation}
The spatial domain is set to
{\color{black}$\Omega=(0,70)\times(0,35)$}, and the boundary conditions
are periodic. The
parameter $\alpha_0$ represents the approach velocity of the
dissipative solitons (and is set to $-0.4$), while $\alpha_5$ is the cubic
cross-coupling parameter (taken equal to $0.5$).
The remaining parameters are set to $\alpha_1=0.125$,
$\beta_1=0.5$, $\alpha_2=-0.9$,
$\alpha_3=1$,
$\beta_3=0.8$, $\alpha_4=-0.1$, and $\beta_4=-0.6$.
The semidiscretization in space is performed with a Fourier pseudospectral
approach with $n_1=700$ modes along the first direction and 
$n_2=350$ along the second one (the total number of degrees of freedom is
then $N=700\cdot350$).

\begin{figure}[!htb]
  \centering
  \input{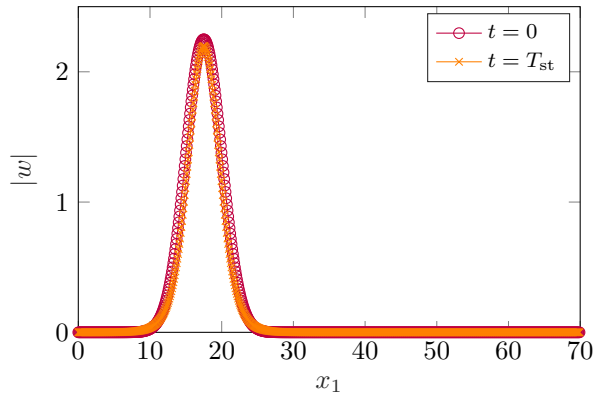}
  \caption{Plot of $\lvert w \rvert$ obtained by simulating the one-dimensional
    cubic-quintic CGL equation~\eqref{eq:cqgl1d}
    with periodic boundary conditions 
  up to final time $T_{\mathrm{st}}=25$. The time integrator is \textsc{if4} 
  with $m=10000$ time steps. The number of spatial discretization
  points is $n_1=700$.
  The wall-clock time of this simulation is 2.62 seconds.
  The obtained state
  is then employed to build the initial condition for the 2D coupled 
  model~\eqref{eq:coupled}.}
  \label{fig:steady_coupled}
\end{figure}

Following the description in Reference~\cite{DB14} (see also the discussion
in References~\cite{ZY15b,M22}), the initial conditions $u_0$ and $v_0$ are 
obtained starting from stationary solitons as follows. 
First, we solved the one-dimensional cubic-quintic CGL equation
\begin{equation}\label{eq:cqgl1d}
  \left\{
  \begin{aligned}
  \partial_tw&=(\alpha_1+\rmi\beta_1)\partial_{x_1x_1} w+
  \alpha_2w+
  (\alpha_3+\rmi\beta_3)\lvert w\rvert^2w+
  (\alpha_4+\rmi\beta_4)\lvert w\rvert^4w,\\
  w_0(x_1)&=\delta\exp\left(-\frac{(x_1-\chi)^2}{2\sigma^2}\right),
  \end{aligned}\right.
\end{equation}
completed with periodic boundary conditions, in the domain $(0,70)$ 
up to $T_{\mathrm{st}}=25$. 
The discretization in space is again performed with a Fourier pseudospectral
method
using $n_1=700$ modes. As time integrator we employ the \textsc{if4}
method with $m=10000$ time steps. The parameters of the initial
condition $w_0$ are set to $\delta=2.25$, $\chi=17.5$, and $\sigma=2.5$.
We selected $T_{\mathrm{st}}=25$ as final simulation time since, after that
time, we did not observe qualitatively variations
of $\lvert w(t,x_1)\rvert$ anymore
(see Figure~\ref{fig:steady_coupled} for the plot of the initial and final
states).
Notice that performing such a simulation is inexpensive in terms of wall
clock time (in fact, it took less than three seconds to obtain the results).
Then, we  set the initial states of the coupled equations~\eqref{eq:coupled}
to $u_0(x_1,x_2)=w(T_{\mathrm{st}},x_1)$ and $v_0(x_1,x_2)=w(T_\mathrm{st},70-x_1)$.
The discretization in space of the coupled PDEs~\eqref{eq:coupled}
yields a system of ODEs~\eqref{eq:ODEs} with
\begin{equation*}
  K=\begin{bmatrix}
  K^{\bu} & 0\\
  0 & K^{\bv}
  \end{bmatrix},
\end{equation*}
where $K^{\bu}$ and $K^{\bv}$ are matrices in Kronecker form
which discretize the linear operators of the equation for the unknown
$u$ and the unknown $v$, respectively. Due to the
block diagonal structure of $K$, the computation of the action
of the operator and its matrix exponential can be independently performed on the
components $\bu$ and $\bv$. This is true also for similar systems with
more than two coupled equations (see, for instance, Reference~\cite{SMS06}).

\begin{table}[!htb]
  \centering
  {\small
  \begin{tabular}{llll|llll}
    \hline
    \multicolumn{4}{l|}{\textsc{if2}} &                      \multicolumn{4}{l}{\textsc{if4}} \\                     
        steps & time (s) & error & order &                  steps & time (s) & error & order \\            
    \hline                                              
    5000 & $1.428\rme{+02}$ & $1.255\rme{-04}$ & ---  &   500 & $3.344\rme{+01}$ & $2.079\rme{-05}$ & ---  \\  
    7000 & $1.982\rme{+02}$ & $6.431\rme{-05}$ & 1.987 &  1000 & $6.712\rme{+01}$ & $1.305\rme{-06}$ & 3.994 \\
    9000 & $2.514\rme{+02}$ & $3.900\rme{-05}$ & 1.990 &  1500 & $1.003\rme{+02}$ & $2.576\rme{-07}$ & 4.001 \\
    11000 & $3.159\rme{+02}$ & $2.614\rme{-05}$ & 1.992 & 2000 & $1.340\rme{+02}$ & $8.130\rme{-08}$ & 4.009 \\
    13000 & $3.759\rme{+02}$ & $1.874\rme{-05}$ & 1.994 & 2500 & $1.677\rme{+02}$ & $3.311\rme{-08}$ & 4.026 \\
    \hline
    \hline
     \multicolumn{4}{l|}{\textsc{strang}}                         &   \multicolumn{4}{l}{\textsc{split4}}\\
     steps & time (s) & error & order                   &   steps & time (s) & error & order\\
                                                          
     5000 & $1.800\rme{+02}$ & $1.008\rme{-04}$ & ---     &   500 & $8.550\rme{+01}$& $4.100\rme{-06}$ & ---     \\
     7000 & $2.536\rme{+02}$ & $5.142\rme{-05}$ & 1.999   &   1000 & $1.712\rme{+02}$ & $2.560\rme{-07}$ & 4.001 \\
     9000 & $3.243\rme{+02}$ & $3.111\rme{-05}$ & 1.999   &   1500 & $2.498\rme{+02}$ & $5.029\rme{-08}$ & 4.013 \\
     11000 & $3.967\rme{+02}$ & $2.083\rme{-05}$ & 1.999  &   2000 & $3.330\rme{+02}$ & $1.569\rme{-08}$ & 4.049 \\
     13000 & $4.678\rme{+02}$ & $1.491\rme{-05}$ & 2.000  &   2500 & $4.158\rme{+02}$ & $6.237\rme{-09}$ & 4.134 \\
    \hline
    \hline
    \multicolumn{4}{l|}{\textsc{rk2}}&
    \multicolumn{4}{l}{\textsc{rk4}}\\
        steps &  & & & 
    steps & time (s) & error & order \\
    \hline
    1062 & $\times$ & $\times$ & ---  &
    562 & $\times$ & $\times$ & --- \\
    1562 & $\times$ & $\times$ & $\times$ &
    1062 & $6.674\rme{+01}$ & $2.574\rme{-03}$ & $\times$\\
    2062 & $5.780\rme{+01}$ & $2.106\rme{-03}$ & $\times$ &
    1562 & $9.886\rme{+01}$ & $5.288\rme{-07}$ & 22.01\\
    2562 & $7.232\rme{+01}$ & $1.341\rme{-03}$ & 2.080  &
    2062 & $1.303\rme{+02}$ & $1.752\rme{-07}$ & 3.978\\
    3062 & $8.601\rme{+01}$ & $9.295\rme{-04}$ & 2.055 &
    2562 & $1.577\rme{+02}$ & $7.393\rme{-08}$ & 3.974\\
    \hline
  \end{tabular}%
  }
  \caption{Number of time steps, wall-clock time (in seconds),
      relative error at final time $T=3$, and observed numerical order of convergence
      for the solution of the 2D coupled cubic-quintic CGL equations~\eqref{eq:coupled} with periodic
      boundary conditions
      and different integrators. The total number of degrees of freedom is
      $N=700\cdot350$. In the table, the cross symbol $\times$ 
      means NaN/Inf value returned.
      See also Figure~\ref{fig:cgl_cubicquintic_coupled_2D_fourier_cpudiag} for a graphical representation.}
  \label{tab:cgl_cubicquintic_coupled_2D_fourier}
\end{table}
\begin{figure}[!htb]
  \centering
%
%
%
\begin{tikzpicture}

\begin{axis}[%
width=2.6in,
height=1.7in,
at={(0.758in,0.481in)},
scale only axis,
xmode=log,
xmin=10,
xmax=1000,
xminorticks=true,
xlabel style={font=\color{white!15!black}},
ymode=log,
ymin=1e-10,
ymax=0.01,
yminorticks=true,
ylabel style={font=\color{white!15!black}},
xlabel = {Wall-clock time (s)},
ylabel = {Relative error},
axis background/.style={fill=white},
legend style={at={(0.02,0.02)}, anchor=south west, legend cell align=left, align=left, draw=white!15!black, font=\footnotesize}
]
\addplot [color=blue, line width=1pt, mark size = 2pt,mark=triangle, mark options={solid, rotate=180, blue}]
  table[row sep=crcr]{%
142.788092	0.000125500337608311\\
198.215791	6.4309233716107e-05\\
251.378925	3.89964423923758e-05\\
315.92459	2.61447826729062e-05\\
375.866098	1.87387473458654e-05\\
};
\addlegendentry{\textsc{if2}}

\addplot [color=cyan, line width=1pt, mark size = 2pt,mark=square, mark options={solid, cyan}]
  table[row sep=crcr]{%
33.436845	2.07894962393836e-05\\
67.117259	1.30505860868947e-06\\
100.286643	2.57639395247736e-07\\
134.004374	8.12999731408415e-08\\
167.685022	3.31110698096309e-08\\
};
\addlegendentry{\textsc{if4}}

\addplot [color=red, line width=1pt, mark size = 2pt,mark=o, mark options={solid, red}]
  table[row sep=crcr]{%
179.958339	0.000100787974953597\\
253.612325	5.14248202142617e-05\\
324.346609	3.11095240842515e-05\\
396.690441000001	2.08256809570219e-05\\
467.693376000001	1.49108650844333e-05\\
};
\addlegendentry{\textsc{strang}}

\addplot [color=magenta, line width=1pt, mark size = 2pt,mark=asterisk, mark options={solid, magenta}]
  table[row sep=crcr]{%
85.503363	4.09967902956841e-06\\
171.217066	2.55976866995751e-07\\
249.840923	5.02911848532878e-08\\
332.976643000001	1.56891289841737e-08\\
415.766348	6.23667506504278e-09\\
};
\addlegendentry{\textsc{split4}}

\addplot [color=black, line width=1pt, mark size = 2pt,mark=triangle, mark options={solid, black}]
  table[row sep=crcr]{%
5.7796e+01	2.1060e-03\\
7.2324e+01	1.3406e-03\\
8.6011e+01	9.2945e-04\\
};
\addlegendentry{\textsc{rk2}}

\addplot [color=green, line width=1pt, mark size = 2pt,mark=diamond, mark options={solid, green}]
  table[row sep=crcr]{%
66.736728	2.574318832852576e-03\\
98.860574	5.287924859396888e-07\\
130.316957	1.752088300148676e-07\\
157.71747	7.393045377583136e-08\\
};
\addlegendentry{\textsc{rk4}}

\end{axis}
\end{tikzpicture}%
  \caption{Results for the simulation of the 2D coupled cubic-quintic CGL
    equations~\eqref{eq:coupled} with periodic
    boundary conditions and
  $N=700\cdot350$ spatial discretization points. The number
  of time steps for each integrator is reported in
  Table~\ref{tab:cgl_cubicquintic_coupled_2D_fourier}. The final simulation time is $T=3$.}
  \label{fig:cgl_cubicquintic_coupled_2D_fourier_cpudiag}
\end{figure}

As did for the numerical examples in the previous sections, we first integrate
the system of cubic-quintic CGL equations~\eqref{eq:coupled} using different
time marching schemes and varying number of time steps.
{\color{black}The nonlinear fluxes in the splitting methods are
approximated by a single step of the \textsc{rk4} scheme.} For the simulations,
the final time is set to $T=3$. The results, together with the specific number
of time steps employed for each method, are reported in
Table~\ref{tab:cgl_cubicquintic_coupled_2D_fourier} and graphically depicted
in Figure~\ref{fig:cgl_cubicquintic_coupled_2D_fourier_cpudiag}.
As we can observe, the most efficient method is also in this case the
fourth-order Lawson scheme \textsc{if4}. In fact, it does not show a
stability restriction
due to the stiffness (as the explicit \textsc{rk4} {\color{black} and 
\textsc{rk2} methods do}) and
overall surpasses all the methods under consideration (in particular
\textsc{split4}, which requires more than twice the wall-clock time needed
by \textsc{if4} to reach about the same level of accuracy).

\begin{figure}[!htb]
  \centering
  \input{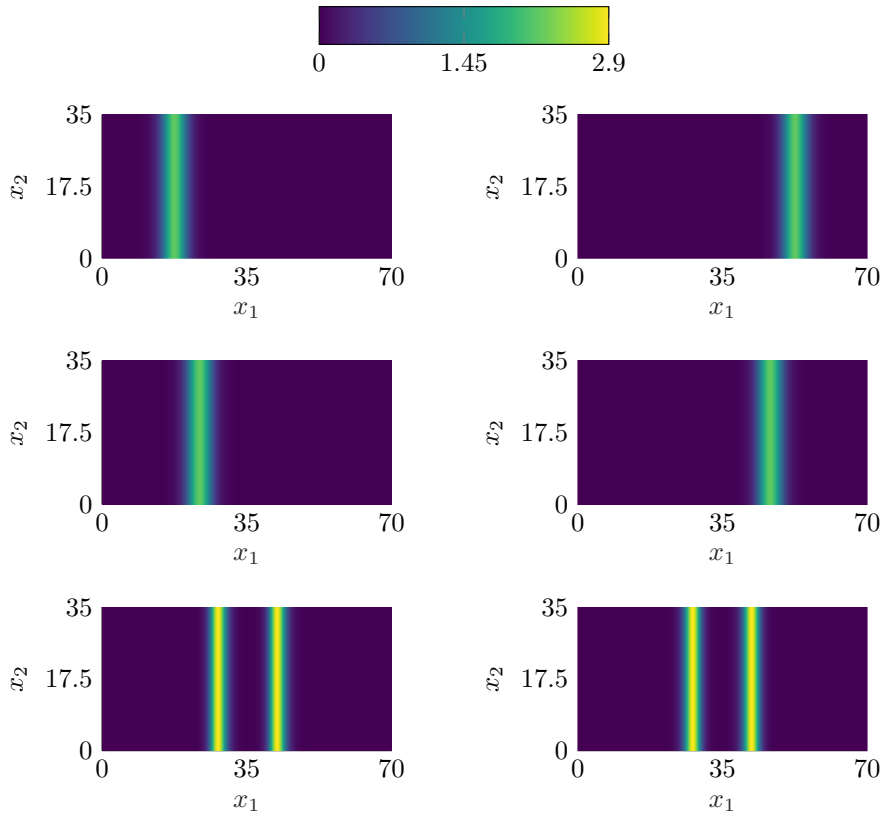}
  \caption{Dynamics of $\lvert u \rvert$ (left) and $\lvert v \rvert$ (right)
    for the 2D
    coupled cubic-quintic CGL equations~\eqref{eq:coupled} with periodic
    boundary conditions and
  $N=700\cdot350$ spatial discretization points. The time integrator is 
  \textsc{if4}. The initial condition is depicted (top), while the 
  simulations are performed with final times $T=15$ (center) and $T=30$ (bottom).
  The
  number of time steps is equal to 
  $m=1000$ and $m=2000$, respectively. The wall-clock times
  of the simulations are $60.89$ and $125.49$ seconds, respectively.}
  \label{fig:cgl_cubicquintic_coupled_2D_fourier_dyn}
\end{figure}

Finally, we perform a last experiment on the coupled model 
by employing \textsc{if4} to show the dynamics of 
$\lvert u \rvert$ and $\lvert v \rvert$
(see Figure~\ref{fig:cgl_cubicquintic_coupled_2D_fourier_dyn}). To this aim,
we consider the final times $T=15$ and $T=30$ with a number of time steps equal
to $m=1000$ and $m=2000$, respectively. As expected, the
advection terms in the equations move the solitons in opposite directions.
When they approach, they collide and start interacting. As thoroughly explained
in Reference~\cite{DB14}, the choice of the parameters in the equations 
should bring to a stationary state constituted by two solitons in the same position
for both the $u$ and the $v$ components. This is indeed verified by the numerical
simulation performed, which after $T=30$ shows qualitatively the same figure.

\section{Conclusions}\label{sec:conc}
We have highlighted, by conducting extensive numerical experiments in a variety
of situations, that complex Ginzburg--Landau equations
{\color{black}defined on Cartesian product domains can be efficiently
  integrated in time by using high-order exponential-type integrators
  with constant time step size.
  We have shown that
  the proposed methods are more efficient (in terms of computational time)
than the explicit fourth-order
   Runge--Kutta scheme  and the
Strang 
splitting procedure.}
Computational efficiency is achieved by suitably exploiting
the underlying Kronecker-sum structure of the problem or
a Fourier 
pseudospectral decomposition in space (when the boundary conditions admit),
both of which 
allow for an effective computation of the involved matrix functions.
We have overall observed that, among the exponential-type schemes, the more the 
nonlinearity becomes involved (i.e., in the cubic-quintic and in the coupled
cases), the better Lawson methods perform compared to splitting-based integrators.
As possible interesting future investigations, we mention
{\color{black}the realization of variable step size integrators,}
the simulation of
further models (the cubic-quintic-septic model proposed in Reference~\cite{DK18},
for instance), and the employment of the proposed techniques in the context
of fractional models (see References~\cite{ZOG21,RDL24}).

\section*{Declaration of competing interest}
The authors declare that they have no known competing financial interests
or personal relationships that could have appeared to influence the
work reported in this paper.
\section*{Acknowledgments}
The authors are members of the Gruppo Nazionale
Calcolo Scientifico-Istituto Nazionale di Alta Matematica (GNCS-INdAM).
Fabio Cassini received financial support from
the Italian Ministry of University and Research (MUR)
with the PRIN Project 2022 No.~2022N9BM3N
``Efficient numerical schemes and optimal control methods for
time-dependent partial differential equations''.
\bibliography{GL24}
\bibliographystyle{elsarticle-num}

\end{document}